\documentstyle[11pt,twoside]{article}

\input{mssymb}

\makeatother

\newcommand{\D}{{\Bbb D}}
\newcommand{\N}{{\Bbb N}}
\newcommand{\R}{{\Bbb R}}

\newcommand{\K}{{\cal K}}
\newcommand{\LL}{{\cal L}}
\newcommand{\CalF}{{\cal F}}
\newcommand{\CalH}{{\cal H}}
\newcommand{\CalK}{{\cal K}}
\newcommand{\CalL}{{\cal L}}
\newcommand{\CalN}{{\cal N}}

\newcommand{\CalU}{{\cal U}}

\newcommand{\Id}{I\mkern-1mud}
\newcommand{\eps}{\varepsilon}
\renewcommand{\epsilon}{\varepsilon}
\newcommand{\bea}{\begin{eqnarray*}}
\newcommand{\eea}{\end{eqnarray*}}
\newcommand{\beq}{\begin{equation}}
\newcommand{\eeq}{\end{equation}}
\newcommand{\rest}[2]{#1\raisebox{-0.5ex}{\mbox{$\mid_{#2}$}}} 
\newcommand{\Lp}{L_p}
\newcommand{\lp}{\ell_p}
\newcommand{\llq}{\ell_q}

\newcommand{\lpcpn}{\ell_p^{\vphantom{n}}(c_p^n)}
\newcommand{\iy}{\infty}
\newcommand{\dopu}{{:}\allowbreak\ }
\newcommand{\teni}{\widehat{\otimes}_{\eps}}

\newcommand{\linq}{\mathop{\overline{{\rm lin}}}\nolimits}
\newcommand{\lin}{\mathop{{\rm lin}}\nolimits}
\newcommand{\ran}{\mathop{\rm ran}\nolimits}
\newcommand{\limp}{\mathop{\rm \vphantom{p}\lim}} 

\newcommand{\mysec}[2]{\section*{\normalsize\hfil\sc {#1}. {#2}\hfill}%
\noindent\setcounter{theo}{0}\setcounter{section}{#1}%
\typeout{#1. #2}}%

\newcommand{\Proof}{\par\noindent{\em Proof. }}
\newcommand{\eop}{\nopagebreak\hspace*{\fill}$\Box$}

\newtheorem{theo}{Theorem}[section]
\newtheorem{lemma}[theo]{Lemma}
\newtheorem{cor}[theo]{Corollary}
\newtheorem{problem}[theo]{Problem}
\newtheorem{prop}[theo]{Proposition}
\newtheorem{df}[theo]{Definition}

\newcounter{abc}   % Counter fr statements-environment wird deklariert
\newcounter{iiiii} % Counter fr aequivalenz-environment wird deklariert

\newenvironment{aequivalenz}
{\setcounter{iiiii}{0}
\begin{list}%
{{\rm (\roman{iiiii})}}%  Falls die items nicht angegeben sind: i)u.s.w.
{\usecounter{iiiii}
\parsep=0pt plus 1pt
\topsep=1pt plus 2pt minus 1pt
\itemsep=1pt plus 2pt minus 1pt
\leftmargin=3\baselineskip
\labelsep=.6\baselineskip
\labelwidth=2.4\baselineskip
\rightmargin 0pt}%
}%               Das war das zweite Argument von "newenvironment"
{\end{list}}

\newenvironment{statements}%
{\setcounter{abc}{0}
\begin{list}%
{{\rm (\alph{abc})}}%  Falls die items nicht angegeben sind: (a) u.s.w.
{\usecounter{abc}
\parsep=0pt plus 1pt
\topsep=1pt plus 2pt minus 1pt
\itemsep=1pt plus 2pt minus 1pt
\leftmargin=3\baselineskip
\labelsep=.6\baselineskip
\labelwidth=2.4\baselineskip
\rightmargin 0pt}%
}%               Das war das zweite Argument von "newenvironment"
{\end{list}}

{\begin{list}%
{$\bullet$}%  Falls die items nicht angegeben sind: Punkte
{\leftmargin=5\baselineskip
\labelsep=\baselineskip
\labelwidth=2.5\baselineskip
\parsep=0pt plus 1pt
\topsep=1pt plus 2pt minus 1pt
\itemsep=1pt plus 2pt minus 1pt
\rightmargin 0pt}%
}%               Das war das zweite Argument von "newenvironment"
{\end{list}}

{\begin{list}%
{$\bullet$}%  Falls die items nicht angegeben sind: Punkte
{\leftmargin 17mm
\labelsep 3mm
\labelwidth 12mm
\rightmargin 0pt}%
}%               Das war das zweite Argument von "newenvironment"
{\end{list}}

\newif\ifrefsc
\let\thebibliographyalt=\thebibliography                                %
\def\thebibliography#1                                                  %
 {\ifrefsc
 \section*{\normalsize\hfil\sc References\hfill}
 \small
 \list                    %
 {[\arabic{enumi}]}{\settowidth\labelwidth{[#1]}\leftmargin\labelwidth  %
 \advance\leftmargin\labelsep                                           %
 \usecounter{enumi}}                                                    %
 \def\newblock{\hskip .11em plus .33em minus .07em}                     %
 \sloppy\clubpenalty4000\widowpenalty4000                               %
 \sfcode`\.=1000\relax                                                  %
 \else\thebibliographyalt{#1}\fi}                                         %

\refsctrue

\begin{document}
\begin{center}
{\Large\bf Property~$(M)$, $M$-ideals, and almost isometric structure
of Banach spaces}
\\[33 pt]
{\sc Nigel J.~Kalton {\it and\/} Dirk Werner}
\end{center}
\bigskip
\begin{quote}
\small
{\sc Abstract.}		
We study $M$-ideals of compact operators by means of the property~$(M)$
introduced in \cite{Kal-M}. Our main result states for a separable
Banach space $X$ that the space of compact operators on $X$ is an
$M$-ideal in the space of bounded operators if (and only if) $X$ does
not contain a copy of $\ell_{1}$, has the metric compact
approximation property, and has property~$(M)$. The investigation of
special versions of property~$(M)$ leads to results on almost isometric
structure of some classes of Banach spaces. For instance, we give a
simple necessary and sufficient condition for a Banach space to embed
almost isometrically into an $\ell_{p}$-sum of finite-dimen\-sional
spaces resp.\ into $c_{0}$, and for $2<p<\iy$ we prove that a
subspace of $L_{p}$ embeds almost isometrically into $\ell_{p}$ if
and only if it does not contain a subspace isomorphic to $\ell_{2}$.
\end{quote}
\bigskip
%
%
%	Section 1
%
%
\mysec{1}{Introduction}%
In 1972 Alfsen and Effros \cite{AlEf} introduced the notion of an
$M$-ideal in a Banach space; a
subspace $J$ of a Banach space $X$ is called an $M$-ideal if there
is an $\ell_1$-direct decomposition of the dual space $X^*$ of $X$
into the annihilator $J^\bot$ of $J$ and some subspace $V\subset X^*$:
$$
X^* = J^\bot \oplus_1 V.
$$
Since that time
$M$-ideal theory has proved useful in Banach
space geometry, approximation theory and harmonic analysis;
see \cite{HWW} for a detailed account.

A number of authors have studied the $M$-ideal structure in $\LL(X,Y)$,
the space of bounded linear operators between Banach spaces $X$ and $Y$,
with special emphasis on the question whether $\K(X,Y)$, the subspace
of compact operators, is an $M$-ideal; see e.g.\ 
\cite{Beh6}, \cite{ChJo1}, \cite{ChJo2}, \cite{GKS}, 
\cite{Kal-M}, \cite{KalW1}, \cite{Lim3}, \cite{OjaCR}, 
\cite{OjaDi}, \cite{PaWe}, \cite{Dirk2}, \cite{WW3}
or Chapter~VI in \cite{HWW}. Examples of $M$-ideals of compact
operators include $\K(\CalH)$ for a Hilbert space $\CalH$ (this is
due to Dixmier \cite{Dix0}), $\K(\lp,\ell_{q})$ for $1<p\le q<\iy$
\cite{Lim3}, and $\K(X)$ for certain renormings of Orlicz sequence
spaces with separable duals \cite{Kal-M}.

Necessary and/or sufficient conditions for $\K(X)$ to be an $M$-ideal
in $\LL(X)$ were studied for instance in \cite{HaLi}, \cite{Kal-M}
and \cite{WW3}. As it turns out, the key notions to tackle this
problem are the properties $(M)$  and $(M^{*})$ introduced in
\cite{Kal-M}. Property~$(M)$ is the
requirement that 
$$
\limsup \|u + x_n\| = \limsup \|v + x_n\| 
$$
whenever $\|u\|=\|v\|$ and $x_n\to0$ weakly, and $(M^*)$ means that
$$
\limsup \|u^* + x_n^*\| = \limsup \|v^* + x_n^*\| 
$$
whenever $\|u^*\|=\|v^*\|$ and $x_n^*\to0$ weak$^*$.
The main result in \cite{Kal-M} shows the equivalence of the
following conditions on a separable Banach space $X$:
(i)~$\K(X)$ is an $M$-ideal in $\LL(X)$; (ii)~$X$ has $(M)$, and
there is a sequence $(K_{n})$ of compact operators such that
$\K_{n}\to\Id$ strongly (that is, $K_{n}x\to x$ for all $x\in X$),
$K_{n}^{*}\to\Id_{X^{*}}$ strongly and $\|\Id-2K_{n}\|\to1$;
(iii)~$X$ has $(M^{*})$, and
there is a sequence $(K_{n})$ of compact operators such that
$\K_{n}\to\Id$ strongly,
$K_{n}^{*}\to\Id_{X^{*}}$ strongly and $\|\Id-2K_{n}\|\to1$.
Here the norm condition on $\Id-2K_{n}$ can be thought of as an
unconditionality property; note that it is fulfilled if $X$ has a
1-unconditional (shrinking) basis.

We now descibe the contents of this paper in detail. 
In section~2 we propose a notion of property~$(M)$ for a single
operator $T\dopu X\to Y$ and point out its relevance for the problem
of giving necessary and sufficient conditions for $\K(X,Y)$ to be
an $M$-ideal in $\LL(X,Y)$. Such a condition is proved for $X=\lp$
in Corollary~\ref{B3}.

Section~3 improves the main result of \cite{Kal-M}, as follows.
It is relatively easy to check that $(M^*)$ implies $(M)$ (see
\cite[Prop.~2.3]{Kal-M}), and it is clear that the converse is 
false (consider $\ell_1$). However in Theorem~\ref{DA} 
we prove that this is 
essentially the only counterexample: If $X$ does not contain $\ell_1$
and has $(M)$, then $X$ has $(M^*)$. It should be noted that 
due to the weak$^*$ compactness of the dual unit ball $(M^*)$
is the technically more convenient property of the two	to work with.
It is also shown in section~3 that the `$\|\Id-2K_n\|\to1$'
condition mentioned above is actually a consequence of 
property~$(M^*)$ and the metric compact approximation property
(as was conjectured in \cite{Kal-M}). Thus we obtain 
for a separable Banach space $X$ that
$\K(X)$ is an $M$-ideal in $\LL(X)$ if and only if
$X$ has $(M)$, does not contain $\ell_1$, and has the
metric compact approximation property.

In section~4 we study a particular version of property~$(M)$, namely
property~$(m_{p})$ which requires that
$$
\limsup \|x+x_{n}\| = (\|x\|^{p} + \limsup\|x_{n}\|^{p})^{1/p}
$$
whenever $x_{n}\to0$ weakly. (The modification for $p=\iy$ is
obvious.) It is known from results in \cite{Kal-M} and \cite{OjaDi}
that, for $1<p\le\iy$, $\K(X\oplus_{p}X)$ is an $M$-ideal in
$\LL(X\oplus_{p}X)$ if and only if $\K(X)$ is an $M$-ideal in $\LL(X)$
and $X$ has $(m_{p})$. Such spaces have been  investigated in
\cite{Kal-M}, \cite{OjaDi} and \cite{PaWe}; and it has been
conjectured (cf.\ \cite[p.~336]{HWW}) that a Banach space with the
above properties must be almost isometric to a subspace of a quotient
of an $\ell_{p}$-sum of finite-dimensional spaces.
(If $\frak Y$ is a class of Banach spaces, we say that $X$ is almost
isometric to a member of $\frak Y$ if for every $\eps>0$ there is a
member $Y_{\eps}$ of $\frak Y$ with Banach-Mazur distance to $X$
less than $1+\eps$.) For separable Banach spaces  a stronger result
than this conjecture is proved in Theorems~\ref{F2} and \ref{F3},
namely, $X$ has $(m_{p})$
if and only if $X$ is almost isometric to a subspace of some $\lp$-sum of
finite-dimensional spaces
if and only if $X$ is almost isometric to a quotient of some $\lp$-sum of
finite-dimensional spaces. As a corollary we obtain an almost
isometric version of a result due to Johnson and Zippin 
\cite{JohZip2} that the
class of subspaces of $\lp$-sums coincides with the class of
quotients of $\lp$-sums. Previously the validity of the above
conjecture has been proved only for special cases if $p=\iy$; see
\cite{WW4} and \cite[p.~337]{HWW}. Theorem~\ref{F5} now presents a
necessary and sufficient condition for a Banach space to embed almost
isometrically into $c_{0}$.

The following two sections deal with property~$(M)$  and in
particular property~$(m_{p})$ for subspaces of $\Lp[0,1]$ (section~5)
resp.\ for subspaces of the Schatten classes $c_{p}$ of operators on
a separable Hilbert space $\CalH$ (section~6). It turns out
(Theorem~\ref{G5}) that a subspace $X$ of $\Lp$, $1<p<\iy$, $p\neq2$,
has $(m_{p})$ if and only if its unit ball is compact for the
topology inherited by $L_{1}$ if and only if $X$ embeds almost
isometrically into $\lp$. Natural examples of this situation
are the Bergman
spaces. If $2<p<\iy$, Theorem~\ref{G5} allows us to deduce a more 
precise variant of a result of Johnson and Odell \cite{JohOde}: $X$
embeds almost isometrically into $\lp$ if and only if $X$ contains no
hilbertian subspace (Corollary~\ref{G6}). In connection with the
Bergman spaces we also consider the Bloch space $\beta_{0}$ and show
that it is almost isometric to a subspace of $c_{0}$ and that
$\K(\beta_{0})$ is an $M$-ideal in $\LL(\beta_{0})$, thus answering a
question raised in \cite{AETD}.

We derive similar results in the context of the Schatten classes
$c_{p}$ in section~6. Here $\lpcpn$ plays the r\^{o}le of $\lp$, and
for $X\subset c_{p}$ we now obtain that $X$
has $(m_{p})$ if and only if $X$ embeds almost
isometrically into $\lpcpn$ (Theorem~\ref{H3}),
and if $p>2$ this happens if and only if
$X$ does not contain a hilbertian subspace (Corollary~\ref{H4}).
This improves a result of Arazy and Lindenstrauss \cite{AraLin}.
A list of open questions concludes the paper.

The results and proofs in this paper are---for the most
part---formulated for the case of real scalars. They extend, however,
to complex scalars and are in fact tacitly used in this context.
Our notation largely follows \cite{LiTz1}; in particular, $B_X$ denotes
the closed unit ball of a Banach space $X$.
The symbol $d(X,Y)$ stands for the Banach-Mazur distance between two
Banach spaces, i.e., $d(X,Y)=\inf \|T\|\,\|T^{-1}\|$ where the
infimum extends over all surjective isomorphisms $T\dopu X\to Y$.
\medskip

\noindent
{\bf Acknowledgement.} This work was  done while the second-named
author  was  visiting the University of Missouri,
Columbia. It is his pleasure to express his gratitude to all
those who made this stay possible.
%\bigskip
%
%
% Section 2
%
%
\mysec{2}{An operator version of property $(M)$}%
In this section we study $M$-ideals of compact operators between
Banach spaces $X$ and $Y$. An in-depth analysis of the case $X=Y$
has been given in \cite{Kal-M}. We first present an operator version of the
main concept of that paper.
%
%
% Def B1
%
%
\begin{df}\label{B1}
\em An operator $T$ with $\|T\|\le 1$ between Banach spaces $X$ and $Y$
is said to have {\em property\/} 
$(M)$ if for all $x\in X$, $y\in Y$, $\|y\|
\le \|x\|$, and weakly null sequences $(x_n)\subset X$,
\beq\label{eqB1}
\limsup \|y+Tx_n\| \le \limsup \|x+x_n\|. 
\eeq
\end{df}

In \cite{Kal-M} 
a separable Banach space was said to have property $(M)$ if the
identity operator has property $(M)$ of Definition~\ref{B1}.
We also mention that Lemma~2.2 of \cite{Kal-M} implies
that every contractive operator in $\LL(X)$ has $(M)$
if and only if the identity operator has $(M)$.

We shall need the following lemma.
%
%
% Lemma B1a
%%
%
\begin{lemma}\label{B1a}
If $T\dopu X\to Y$ is a contractive operator with $(M)$, $(u_n)\subset
X$ and $(v_n)\subset Y$ are relatively compact sequences with
$\|v_n\|\le \|u_n\|$ for all $n$, and if $(x_n)\subset X$ is weakly
null, then
$$
\limsup \|v_n + Tx_n\| \le \limsup \|u_n + x_n\|.
$$
\end{lemma}
\Proof
Otherwise we would have
$$
\limsup \|v_n + Tx_n\| > \limsup \|u_n + x_n\|.
$$
After passing to subsequences we could deal with convergent sequences,
say $\lim u_n = u$ and $\lim v_n=v$, and obtain
$$
\limsup \|v + Tx_n\| > \limsup \|u + x_n\|
$$
contradicting that $\|v\|\le\|u\|$ and that $T$ has $(M)$.
\eop
\bigskip

The next result is a variant of Theorem~2.4 in \cite{Kal-M}, 
and it has a very similar proof.
%
%
% Theo B2
%
%
\begin{theo}\label{B2}
Suppose that the Banach space $X$ admits a sequence of compact operators
$K_n\in \K(X)$ satisfying
\begin{statements}
\item
$K_nx\to x \quad  \forall x\in X$, 
\item
$K_n^*x^*\to x^* \quad  \forall x^*\in X^*$, 
\item
$\|\Id -2K_n\|\to 1$.
\end{statements}
Let $Y$ be a Banach space. Then $\K(X,Y)$ is an $M$-ideal in $\LL(X,Y)$
if and only if every contraction $T\dopu X\to Y$ has $(M)$.
\end{theo}
\Proof
Suppose first that every contraction has property~$(M)$. To show that
$\K(X,Y)$ is an $M$-ideal in $\LL(X,Y)$ we establish the 3-ball
property of \cite[Th.~6.17]{Lim1} (see also \cite[Th.~I.2.2]{HWW}), 
that is, given $S_1, S_2, S_3\in \K(X,Y)$ with
$\|S_i\|\le 1$, $T\in \LL(X,Y)$ with $\|T\|\le 1$, and $\eps>0$,
there is some compact operator $S$ satisfying
$$
\|S_i + T - S\| \le 1+\eps, \quad i=1,2,3.
$$
In fact, we shall show that $S=TK_n$ will work for large $n$.

Since $K_n^* \to \Id_{X^*}$ strongly and $(K_n)$ is uniformly bounded,
$(K_n^*)$ converges uniformly on relatively compact sets and thus
$\|S_i K_m -S_i\| \le \eps/2$, $i=1,2,3$, for large $m$. We also 
assume that $\|\Id - 2K_m\|\le 1 + \eps/2$. Consider a sequence $(x_n)$
in the unit ball of $X$  such that
$$
\limsup_n \|S_1 K_m x_n + (T-TK_n) x_n \| =
\limsup_n \|S_1 K_m  + T-TK_n \|.
$$
Note that $(\Id-K_n)x_n \to 0$ weakly. By Lemma~\ref{B1a} and the fact
that $T$ has $(M)$ we get the inequality
\bea
\limsup_n \|S_1 K_m x_n + T(\Id - K_n) x_n \|
&\le&
\limsup_n \| K_m x_n + (\Id - K_n) x_n \| \\
&\le &
\limsup_n \| K_m  + (\Id - K_n)  \|,
\eea
and the last quantity is dominated by $\|\Id-2K_m\|$ 
(see \cite[p.~153]{Kal-M}).
Therefore, $\limsup_n \|S_1 +T -TK_n\|< 1+\eps$, and
$$
\lim_r \|S_1 +T -TK_{n_r}\|< 1+\eps
$$
for some subsequence $(K_{n_r})$. After switching to further 
subsequences one obtains the same type of inequality for $S_2$ and
$S_3$ which proves the desired 3-ball property.

On the other hand
suppose that $\K(X,Y)$ is an $M$-ideal in $\LL(X,Y)$.
Let $x\in X$, $y\in Y$, $\|y\| \le \|x\|$. Fix a compact operator
$S$ with $\|S\|\le1$ and $Sx=y$. Let $T\dopu X\to Y$ be a contraction
and suppose that $x_n \to 0$ weakly. For each $\eta >0$ there is a 
compact operator $U\dopu X\to Y$ such that $\|Ux - Tx\|\le \eta$
and $\|S+T-U\|\le 1+\eta$ (see \cite[Theorem~3.1, Remark]{Dirk6}). Thus
\bea
\limsup \|y+Tx_n\| 
&=& \limsup \|Sx + Tx_n\| \\
&=& \limsup \|S(x+x_n) + (T-U)x_n\| \\
&=& \limsup \|S(x+x_n) + (T-U)(x+x_n)\| + \eta \\
&\le& (1+\eta) \limsup \|x + x_n\| + \eta.
\eea
Since $\eta$ was arbitrary, this proves that $T$ has property~$(M)$.
\eop
\bigskip

It is easy to check that (\ref{eqB1})
of Definition~\ref{B1} holds if $X$ and $Y$
both have property~$(M)$. The validity of Theorem~\ref{B2} under
these extra assumptions has already
been observed in \cite{Oja4} (cf.\ \cite[Cor.~VI.4.18]{HWW}); 
but we will later encounter situations where
Theorem~\ref{B2} applies, but $Y$ fails $(M)$ 
(cf.\ Propositions~\ref{B4} and~\ref{B5}).

We also mention that the assumptions on $X$ in Theorem~\ref{B2}
enforce $X$ (and $X^*$) to be separable. Nonseparable versions 
(as in \cite{Oja4} and \cite[Chap.~VI.4]{HWW}) are
easily supplied, and we leave it to the reader to spell these out.
Likewise, we don't formulate a corresponding result based on the 
notion of $(M^*)$ of \cite{Kal-M}, which can be done without difficulty.

Note that the assumptions on $X$ above hold in particular if $X$
has a 1-unconditional shrinking basis; e.g., $X=\lp$ for $1<p<\iy$.
In this case there is an embellished version of Theorem~\ref{B2}
which we state next. As usual, $(e_n)$ denotes the unit vector basis
in $\lp$.
%
%
% Cor B3
%
%
\begin{cor}\label{B3}
Suppose $1<p<\iy$. Then $\K(\lp,Y)$ is an $M$-ideal in $\LL(\lp,Y)$
if and only if  for all $y\in Y$  and all operators
$T\dopu \lp\to Y$ the condition
\beq\label{eqB2}
\limsup \|y+Te_n\| \le (\|y\|^p+\|T\|^p)^{1/p} 
\eeq
holds.
\end{cor}
\Proof
The condition is necessary by Theorem~\ref{B2} since
$\limsup \|x+e_n\| = (\|x\|^p+1)^{1/p}$. Next suppose that 
(\ref{eqB2}) holds.
If $(x_n)$ is a weakly null sequence in $\lp$ and if (\ref{eqB1})
fails for some contractive $T\dopu \lp\to Y$, $x\in \lp$ and 
$y\in Y$ with $\|y\|\le\|x\|$, then we may, after passing
to subsequences, assume that
$$
\lim \|y+Tx_n\| > \lim \|x+x_n\|
$$
and that $\gamma=\lim\|x_n\|$	exists. We may also assume without loss
of generality that $\gamma=1$; thus
$$
\lim \|y+Tx_n\| > (\|x\|^p +1)^{1/p}.
$$
On the other hand,
some subsequence $(\xi_n)$ of the sequence $(x_n)$ is almost
isometrically equivalent  to the unit vector basis of $\lp$. Let
$\Phi\dopu \lp \to \overline{{\rm lin}}\{\xi_1, \xi_2, \ldots\}$
be a $(1+\eta)$-isomorphism mapping $e_n$ onto $\xi_n$. Then (\ref{eqB2})
implies that
$$
\limsup \|y+T\xi_n\| = \limsup \|y+(T\Phi)e_n\| \le 
(\|y\|^p+\|T\Phi\|^p)^{1/p},
$$
thus
$$
\limsup \|y+T\xi_n\| \le (\|x\|^p + (1+\eta)^p)^{1/p}
$$
since $\|T\Phi\|\le 1+\eta$, which is a contradiction if $\eta$ is
small enough.
\eop
\bigskip

Let us mention some situations where this corollary applies, thus
providing new examples of $M$-ideals of compact operators.
%
%
% Prop B4
%
%
\begin{prop}\label{B4}
Let $2\le p < \iy$. Then $\K(\lp, \Lp[0,1])$ is an $M$-ideal in
$\LL(\lp,\Lp[0,1])$.
\end{prop}
\Proof
Since the case $p=2$ is well-known, we only have to take care of
the case $2<p<\iy$. Let $T\dopu \lp\to\Lp[0,1]$ be an
operator. We claim that the sequence $(f_n)=(Te_n)$ tends to 0 in 
measure which will prove that 
$$
\limsup \|y+Te_n\| = \limsup (\|y\|^p + \|Te_n\|^p)^{1/p}  
  \le (\|y\|^p + \|T\|^p)^{1/p}  
$$
whenever $y\in \Lp[0,1]$, and hence our assertion by
Corollary~\ref{B3}.

Suppose that $(f_n)$ does not tend to 0 in measure. Then there are
$\eps>0$, $\delta>0$ and Borel sets $A_n$ of measure $\ge \delta$
such that $|f_n| \ge \eps \chi_{A_n}$ a.e.\ for infinitely many $n$.
There is clearly no loss in generality in assuming this to hold
for all $n$. Then we have for some positive constants $c_i$ 
(the existence of $c_1$ below follows from the Khintchin inequalities,
cf.\ \cite[p.~50]{LiTz2}) and all $N\in\N$
\bea
\|T\|\,N^{1/p} &=& \|T\|\, \left\|\sum_{n=1}^N {\pm e_n} \right\|  \\
&\ge& {\rm ave\:} \left\|\sum_{n=1}^N {\pm Te_n} \right\|  \\
&=& {\rm ave\:} \left\|\sum_{n=1}^N {\pm f_n} \right\|  \\
&\ge& c_1 \left\|\biggl(\sum_{n=1}^N |f_n|^2 \biggr)^{1/2} \right\|  \\
&\ge& c_1 \eps 
\left\|\biggl(\sum_{n=1}^N \chi_{A_n} \biggr)^{1/2} \right\|  \\
&\ge& c_2 N^{1/2}
\eea
(where the average extends over all choices of signs), which enforces
$p\le 2$.

Thus the proof of the proposition is completed.
\eop
\bigskip

Note that $\Lp[0,1]$ fails property $(M)$ unless $p=2$.
Proposition~\ref{B4} does not extend to $1<p<2$	because in this case
$\K(\lp,\Lp[0,1])$ is not even proximinal in $\LL(\lp,\Lp[0,1])$
\cite{BangOd}; proximinality is a well-known property of $M$-ideals, 
see e.g.\ \cite[Prop.~II.1.1]{HWW}. On the other hand,
Proposition~\ref{B4} implies the proximinality of
$\K(\lp,\Lp[0,1])$ in $\LL(\lp,\Lp[0,1])$ for $2<p<\iy$ which was
first proved by Bang and Odell \cite{BangOd}.

Clearly the well-known result that $\K(\lp,\ell_q)$ is an $M$-ideal
in $\LL(\lp,\ell_q)$ for $1<p\le q<\iy$ can also be deduced from
Corollary~\ref{B3}.

Next we consider the Schatten class $c_p$ 
of all compact operators on $\ell_2$ with $p$-summable singular values,
and we denote its subspace consisting of those operators $u$ whose
matrix representations with respect to the unit vector basis in
$\ell_2$ are upper triagonal (i.e., $u(e_k)\in {\rm lin}\{e_1, \ldots,
e_k\}$ for all $k$) by $\it UT_p$.
%
%
% Prop B5
%
%
\begin{prop}\label{B5}
~
\begin{statements}
\item
Let $1<p\le 2$. Then $\K(\lp, {\it UT}_p)$ is an $M$-ideal in
$\LL(\lp, {\it UT}_p)$.
\item
Let $2\le p <\iy$. Then $\K(\lp, c_p)$ is an $M$-ideal in
$\LL(\lp, c_p)$.
\item
The space $\K(c_0,\K(\ell_2)) $ is an $M$-ideal in
$\LL(c_0,\K(\ell_2)) $.
\end{statements}
\end{prop}
\Proof
(a) In the case $p=2$ we are dealing with Hilbert spaces, so the result
is clear. Let $1<p<2$. For a fixed $m\in \N$ consider the orthonormal
projection $\pi\dopu \ell_2\to \ell_2$ mapping onto the first $m$
coordinates. We first claim that the operator
$\Phi_p\dopu c_p \oplus_p c_p \to c_p$, $(y,z)\mapsto y\pi + (z-z\pi)$,
is contractive. Indeed, this is trivially true for $p=1$ by the 
triangle inequality, and it holds for $p=2$ since the Hilbert-Schmidt
norm of an operator on $\ell_2$ is just the square sum norm
of the corresponding infinite matrix. Applying the complex 
interpolation method we deduce the same result in the interval
$1<p<2$; note that the complex interpolation space
$[c_1,c_2]_\theta$ for $\theta = 2-\frac{2}{p}$ is $ c_p$ \cite{PieTri}
and that $[c_1 \oplus_1 c_1,	c_2 \oplus_2 c_2]_\theta =
c_p \oplus_p c_p$ \cite[Th.~5.63]{BerghL}.

That $\Phi_p$ is a contraction can be rephrased in the following 
way (denote the matrix entries of an operator $y\in c_p$ by
$y(k,l)$): If, for some $m$, $y(k,l)=0$ for $l>m$ and
$z(k,l)=0$ for $l\le m$, 
then
$$
\|y+z\|_p \le (\|y\|_p^p + \|z\|_p^p)^{1/p}.
$$

Let us now check that (\ref{eqB2}) of Corollary~\ref{B3} holds. Since
the operators with finite matrix representations are dense
in $\it UT_p$, we may assume that $y(k,l)=0$ for $l>m$ for some
$m$. Also, $z_n := Te_n \to 0$ weakly so that $z_n(k,l)\to 0$ for all
$k,l\le m$. Therefore, with $\pi$ as above, $\|z_n\pi\|_p\to 0$
(recall that the $z_n$ are upper triagonal) and thus
\bea
\limsup \|y+Te_n\|_p &=& \limsup \|y+(z_n-z_n\pi)\|_p \\
&\le& \limsup (\|y\|_p^p + \|z_n-z_n\pi\|_p^p)^{1/p} \\
&\le& (\|y\|_p^p + \|T\|_p^p)^{1/p},
\eea
since $(z_n-z_n\pi)(k,l)=0$ for $l\le m$.

(b) Again, we will check (\ref{eqB2}) of Corollary~\ref{B3}, and we may and
shall assume that $y(k,l)=0$ if $k>m$ or $l>m$, for some $m$.
Let $z_n = Te_n$, and let $\pi$ be as above.

If $z\in c_p$, then $\pi z$ is a finite-rank operator mapping
$\ell_2$ into $\ell_2^m$. Since $c_p(\ell_2, \ell_2^m)$ is isomorphic
to $\ell_2$, we can think of the operator $\Psi\dopu \lp \to c_p$
mapping $x\in \lp$ to $\pi\circ Tx$ as an operator from $\lp$
to $\ell_2$, and since $p>2$, $\Psi$ is compact \cite[Prop.~2.c.3]{LiTz1}. 
Consequently,
$\|\pi\circ Te_n\|_p \to 0$. Likewise, $x\mapsto Tx\circ \pi$ is 
compact, and $\|Te_n \circ \pi\|_p \to 0$. Therefore $y+Te_n$
becomes essentially block-diagonal, and we have
\bea
\limsup \|y+Te_n\|_p &=&
\limsup \|\pi y\pi + (\Id-\pi)(Te_n)(\Id-\pi)\|_p \\
&\le& (\|y\|_p^p + \|T\|_p^p)^{1/p}.
\eea

(c) The proof here is similar to that of part~(b).
\eop
\bigskip

The $M$-ideal property of $\K(\lp,V)$ for some subspaces $V$ of $c_p$
has been investigated by P.~Harmand (unpublished). Since by 
Corollary~\ref{B3}	the class of Banach spaces $Y$ for which
$\K(\lp,Y)$ is an $M$-ideal in $\LL(\lp,Y)$ is hereditary, those
results are contained in the above proposition.

The assertion in part~(c) above corresponds to the limiting
case ${p=\iy}$ in~(b). By contrast, we have the following negative
result which was prompt\-ed by a question of T.~S.~S.~Rao.
%
%
% Prop B6
%
%
\begin{prop}\label{B6}
If $\CalH$ is an infinite dimensional Hilbert space, then the three-fold
injective tensor product $\CalH \teni \CalH \teni \CalH$ is not an 
$M$-ideal in its bidual. In fact, 
$\CalH\teni \CalH \teni \CalH = \K(\CalH,\K(\CalH))$ is not
an $M$-ideal in $\LL(\CalH,\K(\CalH))$.
\end{prop}
\Proof
Since the bidual of $\CalH\teni \CalH\teni \CalH$ 
is isometric with $\LL(\CalH,\LL(\CalH))$,
it suffices to prove the latter assertion. It is also enough
to assume that $\CalH$ is separable and thus represented as $\ell_2$.
We shall show that (\ref{eqB2}) of Corollary~\ref{B3} fails.

Let $T\dopu \ell_2 \to \K(\ell_2)$ be the operator that maps the
sequence $(\alpha_n)\in\ell_2$ to the matrix
$$
\left(
\begin{array}{ccccc}
       0 & \alpha_1 & \alpha_2 & \alpha_3 & \ldots \\
\alpha_1	& 0        & 0        &   0       & \ldots \\
\alpha_2 & 0        & 0        &   0      &  \ldots \\
\alpha_3 & 0        & 0        &          &         \\
\vdots   & \vdots   & \vdots   &          & \ddots
\end{array}
\right).
$$
Then $\|T\|=1$, but since the norm of the matrix
$$
\left(
\begin{array}{cc}
1&1\\
1&0
\end{array}
\right)
$$
is $\frac{1}{2} (1+\sqrt{5})$, we have that
$\|e_1 \otimes e_1 + Te_n \| =	\frac{1}{2} (1+\sqrt{5}) > \sqrt{2}$.
\eop
\bigskip

However, for the space $\it UT_\iy$ of upper triangular compact
operators it follows easily from Corollary~\ref{B3} that 
$\K(\CalH, {\it UT}_\iy)$ is an $M$-ideal in 
$\LL(\CalH, {\it UT}_\iy)$.

A more complete picture of those $n$-fold tensor products of
$\lp$-spaces which are $M$-ideals in their biduals will be given in
\cite{Farmer}.
%\bigskip
%
%
% Section 3
%
%
\mysec{3}{Implications of properties $(M)$ and $(M^*)$}%
The first objective of this section is to show that property~$(M)$
implies $(M^*)$ for separable Banach spaces containing no copies
of $\ell_1$; the definitions of properties $(M)$ and $(M^*)$,
introduced in \cite{Kal-M}, have been recalled in the introduction.
Actually, we will also deal with particular versions of $(M)$ and $(M^*)$
that will be studied in more detail in the next sections, but we 
formulate them already here.
%
%
%		Definition Da1
%
%
\begin{df}\label{Da1} \em
Let $1\le p\le\iy$.
\begin{statements}
\item
A separable Banach space is said to have {\em property\/} $(m_p)$ if
$$
\limsup \|x+x_n\| = \|(\|x\|, \limsup \|x_n\|) \|_p
$$
whenever $x_n \to 0$ weakly.
\item
A separable Banach space is said to have {\em property\/} $(m_p^*)$ if
$$
\limsup \|x^*+x_n^*\| = \|(\|x^*\|, \limsup \|x_n^*\|) \|_p
$$
whenever $x_n^* \to 0$ weak$^*$.
\end{statements}
\end{df}

The following lemmas will lead to the first main result, Theorem~\ref{DA}.
%
%
%		Lemma D1
%
%
\begin{lemma}\label{D1}
Suppose $X$ is a separable Banach space and that $(x_n^*)$ is a bounded
sequence in $X^*$.  Let $F$ be a finite-dimensional
subspace of $X^*$. Then for any $\alpha>\liminf d(x_n^*,F)$ there exist
a subsequence $(u_n^*)$ of $(x_n^*)$ and $f^*\in F$ so that
$\|u_n^*-f^*\|\le \alpha$ for all $n$.
\end{lemma}
\Proof Obvious.
\eop
%
%
%		Lemma D2
%
%
\begin{lemma}\label{D2}
Suppose $X$ is a separable Banach space and that
$(x_n^*)$ is a weak$^{\,*}$-null sequence in $X^*$.  Then for any
finite-dimensional subspace $F$ of $X^*$ we have $\liminf d(x_n^*,F)\ge
\frac12\liminf\|x_n^*\|$.
\end{lemma}
\Proof 
Assume $\alpha$ is such that there exist a subsequence
$(u_n^*)$ of $(x_n^*)$ and $f^*\in X^*$ with $\|u_n^*-f^*\|\le \alpha$.
Suppose $\chi$ is any weak$^*$-cluster point of $(u_n^*)$ in $X^{***}$.
Then $\chi\in X^{\perp}$ and $\|\chi-f^*\| \le \alpha$.  Applying the
canonical projection of $X^{***}$ onto $X^*$ we obtain $\|f^*\|\le
\alpha$
and hence $\|u_n^*\|\le 2\alpha$ for all $n$.  The result follows by
Lemma~\ref{D1}.
\eop
%
%
%		Lemma D3
%
%
\begin{lemma}\label{D3}
Suppose $X$ is a separable Banach space containing
no copy of $\ell_1$ and that
$(x_n^*)$ is a weak$^{\,*}$-null sequence.  Then there are a subsequence
$(u_n^*)$ of $(x_n^*)$ and a weakly null sequence $(w_n)$ in $X$ with
$\|w_n\| \le1$ so that $\lim \langle w_n,u_n^*\rangle \ge \frac14 \liminf
\|x_n^*\|$.
\end{lemma}
\Proof
Let $\alpha=\liminf\|x_n^*\|$.  By Lemma~\ref{D2}, we can pass to a
subsequence $(u_n^*)$ so that $\liminf
d(u_n^*,\lin\{u_1^*,\ldots,u_{n-1}^*\})\ge \alpha/2$.  Hence we can find a
sequence $(x_n)$ in $B_X$ so that $\liminf \langle x_n,u_n^*\rangle \ge
\alpha/2$ and $\langle x_n,u_k^*\rangle =0$ for $k\le n$.  By passing
to a subsequence we can suppose that $(x_n)$ is weakly Cauchy
 by Rosenthal's $\ell_1$-theorem \cite[Th.~2.e.5]{LiTz1}.  
Let $w_n=\frac12(x_n-x_{n+1})$, and the
result follows.
\eop
\bigskip

We say that a Banach space has the weak$^{*}$ Kadec-Klee
property if the weak$^*$ and norm topologies agree on its dual
unit sphere. (Occasionally this has been called property~$(**)$.)
%
%
%		Lemma D4
%
%
\begin{lemma}\label{D4}
Assume $X$ is a separable Banach space containing no
subspace isomorphic to $\ell_1$ and with property $(M)$.  Then $X$ has
the weak$^{\,*}$ Kadec-Klee property, 
$X^*$ has no proper norming subspace and $X^*$ is separable.
\end{lemma}
\Proof
Suppose that $(x_n^*)$ converges weak$^*$ to $x^*$ and that
$\lim\|x_n^*\|=\|x^*\|\neq 0$.  We show that $\lim\|x_n^*-x^*\|=0$. 
Suppose not. Then, after passing to some subsequence, we would have
$\lim\|x_n^*-x^*\|=\alpha >0$.  By Lemma~\ref{D3} there are
a weakly null sequence $(w_n)$ in $B_X$ and a subsequence $(u_n^*)$ of
$(x_n^*)$ so that $\liminf\langle w_n,x_n^*-x^*\rangle \ge \alpha/4$.
Hence $\liminf x_n^*(w_n)\ge \alpha/4$.
Suppose $x\in X$ with $\|x\|=1$.  Then
for any $t\ge 0$ we have
$\liminf x_n^*(x+tw_n)\ge  x^*(x) + t\alpha/4$.  By property~$(M)$ we can
pass to a subsequence $(v_n)$ of $(w_n)$ so that there is an absolute
norm
$N$ on
$\R^2$ with $N(1,0)=1$ and $\lim\|u+tv_n\|=N(\|u\|,|t|)$ for any $u\in
X$ and $t\in\R$.  It follows that $N(1,t)-1\ge \alpha|t|/4$.  
But this implies (cf.\ \cite[Lemma~3.6]{Kal-M}) 
that $\ell_1$ embeds into $X$ contrary to our assumption.  
Hence $X$ has the weak$^{*}$ Kadec-Klee property.  We also deduce
 that $X^*$ has no proper norming subspace and must be
separable.
\eop
%
%
%
%		Theo DA
%
%
\begin{theo}\label{DA} 
Let $X$ be a separable Banach
space with property~$(M)$ and containing no copy of $\ell_1$.  
Then $X$ has property~$(M^*)$.
Furthermore if $X$ has property $(m_p)$ where $1<p\le\infty$, then $X$
has property $(m^*_q)$ where $1/p+1/q=1$.
\end{theo}
\Proof
By Lemma~\ref{D4}, $X^*$ is separable.  It therefore suffices, for
the first part, to prove that if $(x_n^*)$
is a weak$^*$-null sequence in
$X^*$ with the property that $\lim \|x^*+x_n^*\|$ exists for all $x^*\in
X^*$, then the function 
$$
\phi(x^*)=\lim\|x^*+x_n^*\|
$$ 
satisfies
$\phi(x^*)=\phi(y^*)$ whenever $\|x^*\|=\|y^*\|$.   For the second part
we will need additionally to show that
$\phi(x^*)=(\|x^*\|^q+\phi(0)^q)^{1/q}$.

We observe first that $\phi$ is a convex and norm-continuous function: 
in fact
$|\phi(x^*)-\phi(y^*)|\le \|x^*-y^*\|$.  We also note $\phi(x^*)\ge
\|x^*\|-\phi(0)$.  For $\tau\ge 0$ we define
$g(\tau)=\inf\{\phi(x^*)\dopu \|x^*\|=\tau\}$.  
Then $g$ is also continuous: in
fact $|g(\tau)-g(\sigma)|\le |\tau-\sigma|$ and $g(\tau)\ge \tau-g(0)$.
It follows that $g$ attains its minimum and that there exists $\tau_0$ so
that $g(\tau_0)\le g(\tau)$ for all $\tau$ with strict inequality if
$\tau_0<\tau$.

Now suppose  that $\tau \ge 0$ and $u^*$ is a weak$^*$-strongly exposed
point of $B_{X^*}$; this means there 
exists $u\in B_X$ with $u^*(u)=1$ such
that  $\lim\|u_n^*-u^*\|=0$ if  $\lim\|u_n^*\|=\lim u_n^*(u)=1$.
The existence of such points follows from the separability of $X^*$, 
see e.g.\ \cite[p.~80]{PheLNM}.
Let $Z=\{x^*\in X^*\dopu  x^*(u)=0\}$.  We will define
$$
\theta=\inf\{\phi(\tau v^*)\dopu  v^*\in u^*+Z\}.
$$  
Let $(F_n)$ be an increasing sequence of
finite-dimensional subspaces of $Z$ so that $\bigcup F_n$ is dense in $Z$.
We claim that $\liminf_{n\to\infty}d(\tau u^*+x_n^*,F_k)\ge \theta$ for
each
$k$.  Indeed if for some $k$ this fails, a simple compactness argument
gives a subsequence $(\xi_n^*)$ of $(x_n^*)$ and $f^*\in F_k$ so that
$\lim
\|f^*+\tau u^* +\xi_n^*\| <\theta$, i.e., $\phi(\tau u^*+f^*)<\theta$.
This would contradict the choice of $\theta$.

Now we may pass to a subsequence $(y_n^*)$ of $(x_n^*)$ so that we
obtain $\liminf d(\tau u^*+y_n^*,F_n)\ge \theta$.  
Hence we can find $y_n\in B_X$ so that
$f^*(y_n)=0$ for $f^*\in F_n$ and 
\beq\label{eqD1}
\liminf(\tau u^*(y_n) + y_n^*(y_n))\ge\theta.
\eeq
Clearly any weak$^*$ cluster point, say $\chi$,  of $(y_n)$ in
$X^{**}$ satisfies $\rest{\chi}{Z}=0$ 
and so $\chi\in \lin\{u\}$.   In particular, the
sequence $(y_n)$ is relatively 
weakly compact, and passing to a subsequence we can
suppose that $(y_n)$ converges weakly to $\alpha u$ for some $\alpha$.
We thus write $y_n=\alpha u + f_n$ where $(f_n)$ is weakly null.

We first use the sequence to estimate $\phi$.  If $\|x\|\le 1$, then by
property~$(M)$, $\limsup \|\alpha x +f_n\|\le \limsup\|\alpha u +f_n\|
\le 1$ and hence for any
$x^*\in X^*$
$$ \phi(x^*) \ge \alpha x^*(x) +\limsup y_n^*(f_n).$$
Let $\beta=\limsup y_n^*(f_n)$;  then
\beq\label{eqDneu} 
\phi(x^*)\ge |\alpha|\|x^*\| + \beta
\eeq
for all $x^*\in X^*$.
Now if $v^*\in u^*+Z$, then $\|v^*\|\ge v^*(u)=1$, and so we obtain
$$ \theta \ge |\alpha|\tau +\beta.$$
On the other hand	by (\ref{eqD1})
$$ \theta\le \liminf(\tau u^*(y_n) + y_n^*(y_n)) \le\alpha\tau + \beta.$$
We conclude that $\theta=\alpha\tau +\beta$ and $\alpha\ge 0$ if
$\tau>0$.

At this point we specialize to the assumption that $\tau>\tau_0$.
Since $\phi(x^*)\ge \alpha\|x^*\| +\beta \ge\beta$ we also have 
by definition of $\tau_0$ that $\beta\le g(\tau_0)$.  
On the other hand,
$$
\theta\ge\inf\{g(\sigma)\dopu\sigma\ge\tau\}>g(\tau_{0}),
$$
and thus $\alpha>0$.  Now choose $v_n^*\in u^*+Z$ so
that $\phi(\tau v_n^*)<\theta +\frac1n$.  Then
$$ \alpha\tau \|v_n^*\| +\beta <\theta  +\frac1n$$
and so $$\|v_n^*\|-1\le \frac1{\alpha\tau n}.$$
Since $v_n^*(u)= u^*(u)=1$ we conclude that
$\lim \|v_n^*\|=1$, and by choice of $u^*$,  $\lim\|v_n^*-u^*\|=0$.
Thus $\phi(\tau u^*)=\theta=\alpha\tau +\beta$.  
By~(\ref{eqDneu}) this implies that
$\phi(\tau u^*) \le \phi(\tau v^*)$ whenever $\|v^*\|=1$ and hence
$\phi(\tau u^*)=g(\tau)$.

Now from the convexity of $\phi$ we have $\phi(\tau x^*) \le g(\tau)$
whenever
$x^*$ is in the convex hull $C$ of the weak$^*$-strongly exposed points of
$B_{X^*}$.  However since $X^*$ is separable $C$ is weak$^*$-dense in
$B_{X^*}$ \cite[p.~80]{PheLNM}.  By the 
weak$^{*}$ Kadec-Klee property (Lemma~\ref{D4})
$C$ is norm-dense in
$B_{X^*}$.  Hence $\phi(\tau x^*)\le g(\tau)$ whenever $\|x^*\|\le 1$.
In particular $\phi(\tau x^*)=g(\tau)$ whenever $\|x^*\|=1$ and
$\tau>\tau_0$.

We can now conclude the argument for the first part. If
$\tau\le
\tau_0$ and
$\|x^*\|=1$ we have by the above argument
$\phi(\tau
x^*) \le g(\sigma)$ for any $\sigma>\tau_0$.  Hence $\phi(\tau x^*)\le
g(\tau_0)$ and consequently $\phi(\tau x^*)=g(\tau_0)=g(\tau)$.  
Thus in either
case $\phi(\tau x^*)=g(\tau)$, and the first part is proved.

For the second part we observe that by what we have already shown
$\phi(u^*)\le \phi(v^*)$ whenever $\|u^*\|\le \|v^*\|$.  Therefore if
$u^*$ is a weak$^*$-strongly 
exposed point of $B_{X^*}$ and $\tau\ge 0$, we can
deduce that in the construction above we have $\phi(\tau
u^*)=\theta=\alpha\tau+\beta$ without restriction on $\tau$.
Hence by~(\ref{eqD1})
$\phi(\tau u^*)\le \|(\tau,\phi(0))\|_q\liminf\|(\alpha,f_n)\|_p$.
Therefore $g(\tau)\le {(\tau^q+g(0)^q)^{1/q}}$.
For the converse direction we 
perform the above construction for $\tau=0$  to produce
$\alpha$, $(f_n)$ so that $\|\alpha u+f_n\|\le 1$ and $g(0)=\lim
y_n^*(f_n)$.  Then $\limsup\|f_n\|\le 1$.  Now for any $\tau>0$ we have
$\limsup \|\gamma u+\delta f_n\| \le \|(\gamma,\delta)\|_p$ and
$g(\tau)\ge \gamma\tau +\delta g(0)$ whenever $\|(\gamma,\delta)\|_p\le
1$ which yields $g(\tau)\ge {(\tau^q+g(0)^q)^{1/q}}$.  The proof is then
complete.
\eop
\bigskip

In the second part of this section we link property~$(M)$ with the
problem of $M$-ideals of compact operators. The next theorem is an
improvement of \cite[Th.~2.4]{Kal-M}. As in \cite{Kal-M} we call
a sequence $(K_n)$ of compact operators a shrinking compact
approximating sequence if $K_n\to \Id_X$ and $K_n^*\to \Id_{X^*}$
strongly.
%
%
%		Theo EA
%
%
\begin{theo}\label{EA}
Let $X$ be a separable Banach space with property
$(M^*)$
and the metric compact approximation property.  Then $\CalK(X)$ is an
$M$-ideal in $\CalL(X)$.
\end{theo}
\Proof By Proposition 2.3 of \cite{Kal-M} $X$ is an $M$-ideal in $X^{**}$,
$X^*$ is separable,
and a result of Godefroy and Saphar \cite[p.~678 and p.~681]{GoSa} 
shows that $X$ has a
shrinking compact approximating sequence $(K_n)$ with $\|K_n\|=1$.  It
remains to show that condition~(6) of Theorem~2.4 of \cite{Kal-M} 
holds, i.e.,
that we can find a shrinking compact approximating sequence 
$(L_n)$ with
$\lim \|\Id-2L_n\|=1$.  It suffices to show that for given $\epsilon>0$ we
can find a shrinking compact approximating sequence $(S_n)$ with $\limsup
\|\Id-2S_n\|<1+\epsilon$.
This will be achieved in a series of lemmas. For notational
convenience we assume $K_{0}=0$.

Let $\CalN_X$ be the collection of norms on $\R^2$ such that there
is a weak$^*$-null sequence $(x_n^*)\subset X^*$ 
satisfying $\|x_n^*\|=1$ and
$N(\alpha,\beta)=\lim_{n\to\infty}\|\alpha u^* +\beta x_n^*\|$ whenever
$\|u^*\|=1$.  It is easy to see that $\CalN_X$ is a compact subset of
the space of all continuous real-valued function on $\R^2$ with the
topology of uniform convergence on compact sets.

Suppose some sequence $(N_j)_{j=1}^{\infty}$ with $N_j \in\CalN_X$
is given. For
each $n$ we define a norm $\Phi$, depending on the sequence $(N_j)$, on
$\R^{n+1}$ by the inductive formula
$$
\Phi(\xi_0,\ldots,\xi_n)=N_n(\Phi(\xi_0,\ldots,
\xi_{n-1}),\xi_n).$$
We can then define a space $\Lambda(N_j)$  or $\Lambda(\Phi)$ as the
completion of the space $c_{00}$ of finitely supported sequences
under the norm
$$ \Phi(\alpha_0,\ldots,\alpha_n,0,\ldots)=\Phi(\alpha_0,\alpha_1,\ldots,
\alpha_n).$$
The space $\Lambda(\Phi)$ is then isomorphic to a subspace of $X^*$.
The collection of such norms $\Phi$ is denoted $\CalF_X$.
%
%
% Lemma E1
%
%
\begin{lemma}\label{E1}
There exist constants $C>0$, $p<\iy$ so that for any
$\Phi\in \CalF_X$, $m\in \N$ and any disjoint $v_1,\ldots,v_m\in
\Lambda(\Phi)$,
 we have 
$$ \biggl(\sum_{j=1}^m \Phi(v_j)^p\biggr)^{1/p}\le C\Phi\biggl(\sum_{j=1}^m
v_j\biggr).$$
\end{lemma}
\Proof
Let $(N'_j)$ be any sequence which is dense in $\CalN_X$ and
form a sequence $(N_j)$ in which each $N'_j$ is repeated
infinitely often; then $(N_j)$ induces a norm $\Psi\in\CalF_X$. Any
$\Lambda(\Phi)$ for
$\Phi\in\CalF_X$ is
(lattice-)finitely representable in $\Lambda(\Psi)$.  Thus it suffices to
show that $\Lambda(\Psi)$ has a lower $p$-estimate for some $p$.  Now
$\Lambda(\Psi)$ is a modular sequence space $h_{(F_j)}$ where
$F_j(t)=N_j(1,t)-1$, for $j\ge 1$ and $F_0(t)=|t|$ (cf.\ \cite{Kal-M}).  
However,
$\Lambda(\Psi)$ is isomorphic to a subspace of $X^*$ and  thus
contains no copy of $\ell_\iy$.  By a result of Woo \cite{Woo}, $h_{(F_k)}$
coincides with a modular
space $h_{(G_k)}$ where the $G_k$ satisfy a uniform $\Delta_2$-condition,
and this in turn implies that $h_{(G_k)}$ has a lower $p$-estimate (by
calculations similar to those in \cite[p.~139f.]{LiTz2}).
\eop
\bigskip

At this point we suppose $0<\delta<\epsilon/8$ and that $\eta_j>0$ are
chosen for $j\ge 1$ so that $\sum_{j\ge 1}\eta_j <\delta$.
%
%
%  Lemma E2
%
%
\begin{lemma}\label{E2}
There exists a subsequence
$(T_n)_{n=0}^{\infty}=(K_{r_n})_{n=0}^{\infty}$ of
$(K_n)$ so that $r_0=0$ and if we put
 $A_n=\{\sum_{j=1}^n\lambda_j(T_j^*x^*-T_{j-1}^*x^*)\dopu 
|\lambda_j|\le1$,
$\|x^*\|\le 1\}$ and 
$B_n=\{T_m^*x^*-T_l^*x^*\dopu \|x^*\|\le 1$, $ m>l\ge n\}$,
then:
\begin{statements}
\item
$\|\sum_{k=1}^n (-1)^kT_k\| < 3$  for $n=1,2,\ldots$~.
\item
For any
 $v^*\in B_n$ there exists $N\in\CalN_X$ so that if $|\alpha|\le 1$ and
 $u^*\in A_{n-1}$, then
$$N(\|u^*\|,|\alpha|\|v^*\|)-\eta_n\le \|u^*+\alpha v^*\|\le
N(\|u^*\|,|\alpha|\|v^*\|)+\eta_n.
$$
\end{statements}
\end{lemma}
\Proof
We will construct $(r_n)$ by induction so that:
\begin{statements}
\item[{\rm (c)}]
$\|\sum_{k=1}^n(-1)^kT_k\|<3$ and $\|\sum_{k=1}^n(-1)^kT_k +
(-1)^{n+1}K_l\|<3$ whenever
$l>r_n$.
\item[{\rm (d)}]
If $D_n=\{K_{m}^*x^*-K_l^*x^*\dopu \|x^*\|\le 1,\ m>l\ge r_n\}$, then
for
any $v^*\in D_n$ there exists $N\in\CalN_X$ so that if $|\alpha|\le 1$
and $u^*\in A_{n-1}$, then $|N(\|u^*\|,|\alpha|\|v^*\|)-\|u^*+\alpha
v^*\|\,| \le \eta_n$.
\end{statements}

To begin the induction simply take $r_0=0$.  Then (c) and (d) are
trivial.

Now suppose $r_0,\ldots,r_{n-1}$ have been chosen.  We show that one
can choose $r_{n}$ large enough so that (c) and (d) hold.  We consider
each case separately.  If one cannot satisfy (c) for all large enough
choices of $r_{n}$ then by induction one can find a sequence $(x_k^*)$
in $X^*$ with $\|x_k^*\|=1$ and increasing sequences of
natural numbers $(a_k)$, $(b_k)$ so that $r_{n-1}<a_1<b_1<a_2<b_2<\cdots$,
for which  if $S=\sum_{j=1}^{n-1}(-1)^jT_j$,
$$ \|S^*x_k^* +(-1)^n (K^*_{b_k}-K^*_{a_k})x_k^*\| \ge 3.$$
By passing to a subsequence we can assume that $(x_k^*)$
converges weak$^*$ to some
$x^*$.  Since $S$ is compact, $(S^*x_k^*)$ converges in norm to $S^*x^*$.
Let $d_k=a_k$ or $b_k$.
Then $(K^*_{d_k}(x_k^*-x^*))$ converges weak$^*$ to $0$ and so
\bea
\limsup \|S^*x_k^* \pm 2K^*_{d_k}(x_k^*-x^*)\| &=& \limsup
\|S^*x^*+2K_{d_k}^*(x_k^*-x^*)\| \\
&=& \limsup \|K_{d_k}^*(S^*x^*+2(x_k^*-x^*))\|\\
&\le& \limsup \|S^*x^* + 2(x_k^*-x^*)\|\\
&\le& \max(\|S\|,2) \limsup \|x^*+(x_k^*-x^*)\|\\
&\le& \max(\|S\|,2).
\eea
 By averaging
$$ \limsup \|S^*x_k^* + (K_{b_k}^*-K_{a_k}^*)(x_k^*-x^*)\| < 3.$$
It remains to observe that $\lim \|(K_{b_k}^*-K_{a_k}^*)x^*\|=0$; then we
arrive at a contradiction and conclude that if $r_{n}$ is chosen large
enough (c) will hold.

We turn now to (d).  If this cannot be satisfied by taking $r_{n}$
large enough, then we can find  a
sequence $(x_k^*)$ with $\|x_k^*\|=1$  and two increasing sequences of
positive integers $(a_k)$, $(b_k)$ so that $r_{n-1}<a_1<b_1<a_2<\cdots$ and
if $v_k^*=(K_{b_k}^*-K_{a_k}^*)x_k^*$ then
for every $N\in\CalN_X$
$$ \max_{|\alpha|\le
1}\max_{u^*\in A_{n-1}}|N(\|u^*\|,|\alpha|\|v_k^*\|)-\|u^*+\alpha
v_k^*\|\,|
>\eta_n.$$ Clearly we have that $\liminf \|v_k^*\|>0$.  We may pass to a
subsequence  so that $(\|v_k^*\|)$ converges to some $\beta>0$ and $\lim
\|x^*+ \alpha v_k^*\|= N(\|x^*\|, |\alpha|\beta)$ for every $x^*\in
X^*$, where $N\in\CalN_X$.  Now pick $|\alpha_k|\le 1$ and $u_k^*\in
A_{n-1}$ so that
$$ |N(\|u_k^*\|,|\alpha_k|\|v_k^*\|)-\|u_k^*+\alpha_kv_k^*\|\,| >\eta_n.$$
By passing to a further subsequence we can suppose $(u_k^*)$ converges in
norm to some $u^*$ (since $A_{n-1}$ is norm compact) and that
$(\alpha_k)$ converges to some $\alpha$.  This gives a contradiction and
so (d) holds for all large enough choices of $r_{n}$.  This enables us
to complete the inductive construction, and so the lemma is
proved.
\eop
%
%
%  Lemma E3
%
%
\begin{lemma}\label{E3}
Suppose $(a_k)_{k\ge 0}$ and $(b_k)_{k\ge 0}$ are two
sequences of nonnegative integers so that $0\le a_0<b_0<a_1<b_1<\cdots$.
Suppose further that $\|x^*\|\le 1$.  Then there exists $\Phi\in
\CalF_X$ so that for any finite sequence
$(\lambda_k)_{0\le k\le n}$, with $|\lambda_k|\le 1$, if we let
$v_k^*=(T_{b_k}^*-T_{a_k}^*)x^*$, we have
$$\Phi(\lambda_0\|v^*_0\|,\ldots,\lambda_n\|v^*_n\|)-\delta\le
\left\|\sum_{k=0}^{n}\lambda_kv_k^*\right\|\le
\Phi(\lambda_0\|v_0^*\|,\ldots,\lambda_n\|v_n^*\|)+ \delta.$$
In particular
\beq\label{eqE1}
 \left\|\sum_{k=1}^n \lambda_kv_k^*\right\| \le 
\left\|\sum_{k=1}^nv_k^*\right\|+2\delta.
\eeq
\end{lemma}
\Proof
It follows from (b) in Lemma~\ref{E2} that we can choose
$N_j\in\CalN_X$ for $j\ge 1$ so that if $u^*\in A_{b_{j-1}}$ and
$|\alpha|\le 1$, then
$$
\bigl|\,N_j(\|u^*\|,|\alpha|\|v_j^*\|) - \|u^*+\alpha v_j^*\|\,\bigr|
<\eta_{a_j}.
$$
Choose $\Phi$ corresponding to the sequence
$(N_j)$.  Then we prove by induction that
\bea
\Phi(\lambda_0\|v_0^*\|,\ldots,\lambda_k\|v_k^*\|)-\sum_{j=1}^k\eta_{a_j}
&\le&
\left\|\sum_{j=0}^k\lambda_jv_j^*\right\|     \\
&\le&
\Phi(\lambda_0\|v_0^*\|,\ldots,\lambda_k\|v_k^*\|)+\sum_{j=1}^k\eta_{a_j}
\eea
for $ 0\le k\le n$.  This is trivial for $k=0$.  Assume it is true
for $k=s-1$.  Then
$$ \left|\,\left\|\sum_{j=0}^{s}\lambda_j v_j^*\right\| -
N_s\left(\left\|\sum_{j=0}^{s-1}\lambda_jv_j^*\right\|,
\lambda_{s}\|v_{s}^*\|\right) \right| \le
\eta_{a_s}.$$
The inductive hypothesis follows easily.  The last statement is trivial
from the unconditionality of the norm $\Phi$.
\eop
\bigskip

Now let $V_n=T_{n+1}-T_{n}$.
%
%
%	Lemma E4
%
%
\begin{lemma}\label{E4}
For any finite sequence $(\lambda_k)_{k=0}^n$ we
have 
$$\left\|\sum_{k=0}^n\lambda_kV_k\right\| \le 8\max |\lambda_k|.$$
\end{lemma}
\Proof
For convenience we suppose $n$ is even, say $n=2m;$ we also
suppose  $\max |\lambda_k|=1$.  Let
$\|x^*\|=1$. Then by Lemma~\ref{E3},
$$
\left\|\sum_{k=0}^{m}\lambda_{2k}V^*_{2k}x^*\right\|\le
\left\|\sum_{k=0}^mV_{2k}^*x^*\right\|+2\delta.$$
However by (a) in Lemma~\ref{E2} we have 
$\|\sum_{k=0}^mV^*_{2k}\| <3$.  Hence
$$ \left\|\sum_{k=0}^m\lambda_{2k}V^*_{2k}x^*\right\| \le 4.$$
 Similarly
$$ \left\|\sum_{k=1}^m\lambda_{2k-1}V^*_{2k-1}x^*\right\|\le 4$$ 
and the lemma is
proved.
\eop
\bigskip

It follows now that for any $x^*\in X^*$ the series $\sum_{n=0}^{\infty}
\lambda_nV_n^*x^*$ converges weak$^*$ in $X^*$.  In fact since $X^*$ is
separable the series must converge in norm, although we do not need this
observation.
%
%
%	Lemma E5
%
%
\begin{lemma}\label{E5}
Suppose $C$ and $p$ are given by Lemma~\ref{E1}, and that $t\in\N$.  
Then for $\|x^*\|=1$ and any  sequence $(\lambda_k)$ with
$\max|\lambda_k|\le 1$ we have:
$$ \left\|\sum_{k=0}^{\infty}\lambda_kV_k^*x^*\right\| 
\le 1 + 18Ct^{-1/p} +8\delta + 16t\Delta$$ 
where
$\Delta=\sup_{n\ge 0}|\lambda_n-\lambda_{n+1}|$.
\end{lemma}
\Proof We begin by noting that according to Lemma~\ref{E3}
there exists $\Phi\in\CalF_X$ so
that if $(\alpha_k)_{k=0}^n$ is a finite sequence satisfying
$\max_{0\le k\le n} |\alpha_k|\le 1$, then
\bea
 \Phi(\alpha_0\|V_0^*x^*\|,\ldots,\alpha_n\|V_{2n}^*x^*\|)-\delta 
&\le&
\left\|\sum_{j=0}^n\alpha_jV_{2j}^*x^*\right\|    \\
&\le&
 \Phi(\alpha_0\|V_0^*x^*\|,\ldots,\alpha_n\|V_{2n}^*x^*\|)+\delta.
\eea
We can apply Lemma~\ref{E1} to obtain
for any $n\in\N$
$$ \left(\sum_{s=0}^{t-1}\left(\left\|\sum_{k=0}^n V_{2tk+2s}^*x^*\right\|
-\delta\right)^p\right)^{1/p}
\le C
\left(\left\|\sum_{k=0}^{nt+t-1}V_{2k}^*x^*\right\|+\delta\right)\le 9C.$$
It follows that there exists $0\le s\le t-1$ so that
$$ \left\|\sum_{k=0}^n V_{2tk+2s}^*x^*\right\| \le \delta + 9Ct^{-1/p}.$$
Now, by Lemma~\ref{E3},
$$
\left\|\sum_{k=0}^n \lambda_{2tk+2s}V_{2tk+2s}^*x^*\right\| \le
\left\|\sum_{k=0}^n V_{2tk+2s}^*x^*\right\|+2\delta
\le 9Ct^{-1/p}+3\delta.
$$
For $0\le k\le n+1$ we set
$$ v_k^* =\sum_{2(k-1)t+2s<j<2kt+2s}V_j^*x^*$$
and
$$ w_k^* =\sum_{2(k-1)t+2s<j<2kt+2s}\lambda_jV_j^*x^*,$$
where each range of summation consists only of $j$ with $0\le j\le
2nt+2t-1$.  Then  for $0\le k\le n$ we have
$$ w_k^*-\lambda_{2kt+2s}v_k^* =
\sum_{2(k-1)t+2s<j<2kt+2s}(\lambda_j-\lambda_{2kt+2s})V_j^*x^*.$$
Similarly
$$
w_{n+1}^*-\lambda_{2nt+2s}v_{n+1}^*=\sum_{2nt+2s<j<2nt+2t-1}
(\lambda_j-\lambda_{2nt+2s})V_j^*x^*.$$
Combining we have that
$$
\sum_{k=0}^{n+1}w_k^*-\sum_{k=0}^{n+1}\alpha_kv_k^*=2t\Delta\sum_{j=0}^
{2nt+2t-1}\beta_jV_j^*x^*$$
where $\max_{0\le k\le n+1}|\alpha_k|\le 1$ and $\max_{0\le j\le
2nt+2t-1}|\beta_j|\le 1$.
Thus by (\ref{eqE1}) and Lemma~\ref{E4}
$$ \left\|\sum_{k=0}^{n+1}w_k^*\right\| 
\le \left\|\sum_{k=0}^{n+1}v_k^*\right\| +2\delta
+16t\Delta.$$

Now
\bea
\left\|\sum_{j=0}^{2nt+2t-1}\lambda_jV_j^*x^*\right\| 
&\le& 
\left\|\sum_{k=0}^{n+1}w_k^*\right\|
+ \left\|\sum_{k=0}^n\lambda_{2tk+2s}V_{2tk+2s}^*x^*\right\| \\
&\le& 
\left\|\sum_{k=0}^{n+1}v_k^*\right\| + 9Ct^{-1/p} +16t\Delta+5\delta\\
&\le& 
\left\|\sum_{j=0}^{2nt+2t-1}V_j^*x^*\right\| 
  + 18Ct^{-1/p} +16t\Delta +6\delta  \\
&\le&
1+ 18Ct^{-1/p}+16t\Delta+8\delta.
\eea
Finally letting $n\to\infty$ we obtain the result.
\eop
\bigskip

We complete the proof of Theorem~\ref{EA}
by showing that if $S_n=\frac1n(T_1+\cdots+T_n)$,
then $\limsup \|\Id-2S_n\| < 1+\epsilon$.
To see this we write
$$ \Id-2S_n = \sum_{k=1}^n \frac{2k-n-2}{n}V_k + \sum_{k=n+1}^{\infty}V_k$$
so that by Lemma~\ref{E5} for any $t\in \N$
$$ \|\Id-2S_n\| \le 1+18Ct^{-1/p} + 32tn^{-1} +8\delta.$$
Hence $\limsup \|\Id-2S_n\| \le 1+ 18Ct^{-1/p}+8\delta$.  However, as $t$
is arbitrary we obtain $\limsup\|\Id-2S_n\| < 1+\epsilon$.
\eop
\bigskip

Combining Theorems \ref{DA} and \ref{EA} we obtain the following
result.
(For the necessity of the conditions below see \cite{HaLi} and 
\cite{Kal-M}.)
%
%
% Corollary Dtheo
%
%
\begin{theo}\label{Dtheo}
Let $X$ be a separable Banach space. Then $\K(X)$ is an $M$-ideal in
$\LL(X)$ if and only if $X$ does not contain a copy of $\ell_1$,
$X$ has $(M)$ and $X$ has the metric compact approximation
property.
\end{theo}

Since a reflexive space with the compact approximation property
automatically has the metric compact approximation property 
\cite{ChJo1} (see also \cite{GoSa}) we can
refine Theorem~\ref{Dtheo} in the case of reflexive $X$.
\begin{cor}\label{reflcor}
Let $X$ be a separable reflexive Banach space. 
Then $\K(X)$ is an $M$-ideal in
$\LL(X)$ if and only if $X$ has $(M)$ and the 
compact approximation property.
\end{cor}
%
%
%  end of part 1
%
%

%
%
% Section 4
%
%
\mysec{4}{Implications of property $(m_p)$}%
In this section we prove that a separable Banach space that does
not contain a copy of $\ell_1$ and
enjoys property $(m_p)$, which was introduced in
Definition~\ref{Da1}, 
is almost isometric to a subspace of
an $\ell_p$-sum of finite-dimensional spaces if $1<p<\iy$ 
(Theorem~\ref{F3}), respectively almost isometric to a subspace
of $c_0$ if $p=\iy$ (Theorem~\ref{F5}).
At the end of this section we point out the relevance
of the $(m_p)$-condition in the theory of $M$-ideals.

We first recall that, 
if $1\le p<\infty$, the space $C_p$ is defined to be  $\ell_p(E_n)$ 
where $(E_n)$ is a sequence of finite-dimensional Banach spaces
dense
in all finite-dimensional Banach spaces in the sense of Banach-Mazur
distance.  
If $p=\infty$ we define $C_{\infty}=c_0(E_n)$.  Now if $(F_n)$
is any sequence of finite-dimensional Banach spaces and $\epsilon>0$, it
is easy to see that there is a 1-complemented subspace $X_0$ of $C_p$ so
that $d(\ell_p(F_n),X_0)<1+\epsilon$, or $d(c_0(F_n),X_0)<1+\epsilon$
when $p=\infty$.
For short we say that $\ell_{p}(F_{n})$  is almost isometric to a
1-complemented subspace of $C_{p}$.
It is clear that
$C_\iy$ is almost isometric to a subspace of $c_0$. 
 The spaces $C_p$ have been investigated by Johnson
and Zippin (\cite{John2}, \cite{John0}, \cite{JohZip}, \cite{JohZip2}).

Although the definition of $C_p$ depends on the choice of the 
sequence $(E_n)$, it is known that any two such choices yield
essentially the same space $C_p$. More precisely, if $C_p$ and $C_p'$
are built on two dense sequences $(E_n)$ and $(E_n')$, then $d(C_p,
C_p')=1$ (see \cite{John2}). 
Since the statements we are going to make are almost 
isometric in character rather than isometric, this ambiguity 
does not affect the following theorems. So we fix one model of $C_p$
once and for all.

Some of the ideas in the proof of the following Theorem~\ref{F2}
are borrowed from those of Johnson and Zippin \cite{JohZip2}.
We also need a lemma.
%
%
%  Lemma F1
%
%
\begin{lemma}\label{F1}
Let $1<p\le\iy$ and let $X$ be a separable Banach space with $(m_p)$ that
does not contain a copy of $\ell_1$.
Let $q$ be the exponent conjugate to $p$, i.e., $1/p + 1/q =1$.
\begin{statements}
\item
If $E$ is a finite-dimensional subspace of $X$ and $\eta>0$, 
then there is a finite-codimensional subspace $U$ of $X$ such that
$$
(1-\eta) (\|x\|^p + \|y\|^p)^{1/p} \le \|x+y\| \le
(1+\eta) (\|x\|^p + \|y\|^p)^{1/p}
$$
for all $x\in E$, $y\in U$.
\item
If $F$ is a finite-dimensional subspace of $X^*$ and $\eta>0$, 
then there is a finite-codimensional weak$^{\,*}$-closed
subspace $V$ of $X^*$ such that
$$
(1-\eta) (\|x^*\|^q + \|y^*\|^q)^{1/q} \le \|x^*+y^*\| \le
(1+\eta) (\|x^*\|^q + \|y^*\|^q)^{1/q}
$$
for all $x^*\in F$, $y^*\in V$.
\end{statements}
\end{lemma}
\Proof
(a)  A compactness argument shows that it is enough to give the proof
in the case where $\dim E=1$, say $E=\lin\{x\}$. 
By Lemma~\ref{D4} $X^*$ is separable; let 
$(u_n^*)$ be a dense sequence in $X^*$. 
Suppose now that the lemma does not hold.
Then there exists some $\eta>0$ so that
there is a bounded sequence  $(y_n)$, $y_n \in U_n = \{x\dopu
u_k^*(x)=0$, for $k=1,\ldots,n\}$, satisfying
$$
\bigl|\, \|x+y_n\| - (\|x\|^p+\|y_n\|^p)^{1/p}\, \bigr| >  \eta.
$$
Since $y_n\to 0$ weakly, this is a contradiction to property $(m_p)$.

The proof of (b) is similar; note that $X$ has $(m_q^*)$ by 
Theorem~\ref{DA}.
\eop
%
%
%		Theo F2
%
%
\begin{theo}\label{F2}
Suppose $1<p\le \infty$, and suppose $X$ is a separable
Banach space not containing a copy of $\ell_1$.  
Then $X$ has property $(m_p)$ if and only if given any
$\epsilon>0$ there is a quotient $X_0$ of $C_p$ so that
$d(X,X_0)<1+\epsilon$.
\end{theo}
\Proof
It is clear by Theorem~\ref{DA}
that -- in the case $p<\iy$ -- 
a space which is almost isometric to a quotient of
$C_p$ has $(m_p)$, since it is reflexive and its dual has $(m_{q})$.
(Here and in what follows we suppose   $1/p + 1/q=1$.)
In the case $p=\iy$  we use a result due to Alspach \cite{Als} to
conclude that a quotient of $C_{\iy}$ is almost isometric to a
subspace of $c_{0}$ and thus has $(m_{\iy})$.

We now prove the converse.
Suppose $\delta>0$ is chosen so that
$\delta<\frac14\epsilon$ and
$\frac{1-\delta}{1+\delta}>1-\frac18\epsilon$.  We will also need a
positive integer $t$ selected so that
$32t^{-1/q}<\frac18\epsilon$.  Suppose
$(\eta_n)$ is a sequence of real numbers satisfying
$0<\eta_n<\frac12\delta$, $\prod_{n\ge 1}(1-\eta_n)>1-\delta$ 
and $\prod_{n\ge 1}(1+\eta_n)<1+\delta$.

Let $(u_n)$ be a dense sequence in $X$.
For $F\subset X^*$ we let $F_\bot =\{x\in X\dopu x^*(x)=0\ \forall
x^*\in F\}$, and we denote the span of $u_1,\ldots,
u_n$ by  $[u_1,\ldots,u_n]$.
  We will inductively choose two
sequences of
subspaces $(F_n)$ and $(F_n')$ of $X^*$ and
subspaces $E(m,n)$, $1\le m\le n$, in $X$ so that:
\begin{statements}
\item
 $\dim F_n<\infty,\ \dim E(m,n)<\infty$ for all $m\le n$.
\item
 $F_n'\subset [u_1,\ldots,u_n]^{\perp}\cap\bigcap_{j\le
k<n}E(j,k)^{\perp}$ is weak$^*$-closed and
$X^{*}=F_1\oplus\cdots\oplus F_n\oplus F_n'$.
\item
 $F'_n = F_{n+1}\oplus F'_{n+1}$.
\item
 If $x^*\in F_1+\cdots+ F_n$ and $y^*\in F'_{n+1}$, then
$$(1-\eta_n)(\|x^*\|^{q}+\|y^*\|^q)^{1/q}  \le
\|x^*+y^*\|\le
(1+\eta_n)(\|x^*\|^{q}+\|y^*\|^q)^{1/q}.$$
\item
 If $x\in (F_1+\cdots +F_n)_{\perp}$ and $y\in \sum_{j\le k<n}E(j,k)$,
then
$$ (1-\eta_n)(\|x\|^p+\|y\|^p)^{1/p}
\le \|x+y\| \le
(1+\eta_n)(\|x\|^p+\|y\|^p)^{1/p}
$$
when $p<\infty$ or
$$ (1-\eta_n)\max(\|x\|,\|y\|) \le \|x+y\| \le
(1+\eta_n)\max(\|x\|,\|y\|)$$
when $p=\infty$.
\item
 We have $(F_1+\cdots+F_{m-1}+F'_n)_{\perp}\subset E(m,n)$ 
 and $E(m,n)\subset (F_1+\cdots +F_{m-2})_{\perp}$ if $1\le m\le
n$.
\item
 If $x^*\in F_m+\cdots+F_n$, then there exists $x\in E(m,n)$ so that
$\|x\|\le 1$ and $x^*(x)\ge (1-\delta)\|x^*\|$.
\end{statements}
(Here we understand $F_1+\cdots +F_r = \{0\}$ if $r<1$.)

Let us describe the inductive construction.  We pick $F_1'=[u_1]^{\perp}$
and $F_1$ to be any complement of $F_1'$ to start the procedure.

Now suppose we have defined $(F_j)_{j\le n}$,  $(F_j')_{j\le n}$ and
$E(j,k)$ for $1\le j\le k\le n-1$ (at the first step no $E$-subspaces
exist).

We first must define the spaces $E(m,n)$ for $1\le m\le n$.  To do this
we notice that if $m\ge 3$, $x^*\in F_1+\cdots +F_{m-2}$ and $y^*\in
F_m+\cdots +F_n$ we have
$$ \|x^*+y^*\| \ge (1-\eta_{m-2})(\|x^*\|^q+\|y^*\|^q)^{1/q}\ge
\Bigl(1-\frac{\delta}2\Bigr)\|y^*\|.$$
Hence we can find a finite-dimensional subspace $G=G(m,n)$ so that
$G\subset (F_1+\cdots +F_{m-2})_{\perp}$  and
$$ \sup_{g\in G,\|g\|=1} y^*(g) \ge (1-\delta)\|y^*\|$$
if $y^*\in F_m+\cdots+F_n$.  We then augment $G$ to $E(m,n)=G +
(F_1+\cdots+ F_{m-1}+F_n')_{\perp}$. For $m=1$ or $m=2$ the procedure 
is similar.

We now turn to the definitions of $F_{n+1}$ and $F'_{n+1}$.   Using 
Lemma~\ref{F1} we see that there
exist finite-dimensional subspaces $H\subset X$ and $K\subset X^*$
so that:
\begin{statements}
\item[(h)]
 If $x\in \sum_{1\le m\le n}E(m,n)$ and $y\in K_{\perp}$, then
$$ (1-\eta_{n+1})(\|x\|^p+\|y\|^p)^{1/p} \le \|x+y\| \le
(1+\eta_{n+1})(\|x\|^p+\|y\|^p)^{1/p}.$$
\item[(i)]
 If $x^*\in F_1+\cdots +F_n$ and $y^*\in H^{\perp}$, then
$$ (1-\eta_n)(\|x^*\|^q +\|y^*\|^q)^{1/q} \le \|x^*+y^*\| \le
(1+\eta_n)(\|x^*\|^q +\|y^*\|^q)^{1/q}.$$
\end{statements}

Let $D$ be any weak$^*$-closed complement of $F_1+\cdots+F_n+K$ which is
contained in $F_n'$ (which is a complement of $F_1+\cdots+F_n$).  Then
let $F'_{n+1}$ be defined as the intersection
$$ F'_{n+1} = D\cap \biggl(H+\sum_{1\le m\le
n}E(m,n)+[u_1,\ldots,u_n]\biggr)^{\perp}.$$
Now there is a subspace $K'$ of $F_n'$ so that $F_1+\cdots
+F_n+K=F_1+\cdots +F_n+K'$.  Pick $F_{n+1}\supset K'$ so that
$F_n'=F_{n+1}\oplus F'_{n+1}$.  This completes the inductive
construction.

Now suppose $(a_n)_{n\ge 0}$ and $(b_n)_{n\ge 0}$ are two increasing
sequences of integers with $a_0\le b_0<a_1\le b_1<a_2\le\cdots$.  First
suppose $b_n+2\le a_{n+1}$ for all $n$.  We observe
 that if
$x_n^*
\in F_{a_n}+\cdots+F_{b_n} $ for all $n$, then by induction, using (d), we
have
$$ \prod_{k=1}^{n-1}(1-\eta_{b_k})\biggl(\sum_{k=1}^n\|x^*_k\|^q
\biggr)^{1/q}
\le \left\|\sum_{k=1}^nx^*_k\right\| \le
 \prod_{k=1}^{n-1}(1+\eta_{b_k})\biggl(\sum_{k=1}^n\|x^*_k\|^q
\biggr)^{1/q}$$
and hence
$$
(1-\delta)\biggl(\sum_{k=1}^n\|x^*_k\|^q\biggr)^{1/q}
\le \left\|\sum_{k=1}^nx^*_k\right\| \le
(1+\delta)\biggl(\sum_{k=1}^n\|x^*_k\|^q\biggr)^{1/q}.$$
Similarly if $b_n+3\le a_{n+1}$ and $x_n\in E(a_n,b_n)$ then by (e) and (f)
and induction
$$
(1-\delta)\biggl(\sum_{k=1}^n\|x_k\|^p\biggr)^{1/p}
\le \left\|\sum_{k=1}^nx_k\right\| \le
(1+\delta)\biggl(\sum_{k=1}^n\|x_k\|^p\biggr)^{1/p}$$
when $p<\infty$ or
$$ (1-\delta)\max(\|x_1\|,\ldots,\|x_n\|)
\le \left\|\sum_{k=1}^nx_k\right\| \le
 (1+\delta)\max(\|x_1\|,\ldots,\|x_n\|)$$
when $p=\infty$.

We also notice that for each $n$ there is a weak$^*$-continuous
projection, $P^*_n$ say, of $X^*$ onto $F_n$ with kernel $F_1+\cdots+
F_{n-1}+F_n'$.  Then $P^*_nP^*_m=0$ if $m\neq n$, and $\ran P_k\subset
E(m,n)$ whenever $m\le k\le n$.  If $x^*\in X$ and $1\le n<\infty$ then
we can write $x^*=\sum_{k=1}^nP_k^*x^*+w^*$ with
$w^*\in F'_n$.   Then
$x^*(u_n)=\sum_{k=1}^nx^*(P_ku_n)$ and $P_ku_n=0$ for $k>n$.

Now if $0\le s\le t-1$ we will define two operators, $T_s$ and $R_s$.
Suppose $p<\infty$.
We define 
$$
T_s\dopu Y_{s}{:=}
\ell_p\bigl(E(4(n-1)t+4s+4,4nt+4s+1)_{n\ge 0}\bigr)\to X
$$
by
$T_s\bigl((x_n)\bigr)=\sum_{n\ge 0}x_n$ 
where the first space with $n=0$ is to be
understood as $E(1,4s+1)$.   We also define
$$
R_s\dopu Z_{s}{:=}\ell_p\bigl(E(4nt+4s+2,4nt+4s+3)
_{n\ge 0}\bigr)\to X
$$ 
by
$ R_s\bigl((x_n)\bigr) =\sum_{n\ge 0}x_n$.
In the case $p=\infty$ we let $Y_{s}$ and $Z_{s}$ be
$c_0$-sums.

It is clear that $(1-\delta)\|\xi\|\le \|T_s\xi\|\le (1+\delta)\|\xi\|$
for $\xi\in Y_{s}$ and similarly $(1-\delta)\|\xi\|\le \|R_s\xi\|\le
(1+\delta)\|\xi\|$ for $\xi \in Z_{s}$.

At the same time we introduce two operators
$T\dopu Y{:=}\ell_p\bigl((Y_s)_{s=0}^{t-1}\bigr)\to X$ 
and $R\dopu Z{:=}\ell_p\bigl((Z_s)_{s=0}^{t-1}\bigr)\to	X$ defined by
$$T(\xi_0,\ldots,\xi_{t-1})=\frac1{t^{1/q}}\sum_{s=0}^{t-1}T_s\xi_s$$
and
$$R(\xi_0,\ldots,\xi_{t-1})=\sum_{s=0}^{t-1}R_s\xi_s.$$
By construction we have 
$$
\|T(\xi_0,\ldots,\xi_{t-1})\| \le 
\frac1{t^{1/q}}\sum_{s=0}^{t-1}\|T_s\|\,\|\xi_s\| \le
(1+\delta)\biggl(\sum_{s=0}^{t-1} \|\xi_s\|^p \biggr)^{1/p}
$$
and thus $\|T\|\le 1+\delta$.  Also, we can consider
$R$ as an operator on the space
$\ell_p\bigl(E(4n+2,4n+3)_{n\ge 0}\bigr)$ given  by
$R\bigl((x_n)\bigr)=\sum_{n\ge 0}x_n$ and deduce that $\|R\|\le 1+\delta$.

Now suppose $x^*\in X^*$ and $R_s^{*}x^*=0$.  Then it is clear that
$P_k^*x^*=0$ whenever $k\equiv 4s+2$ or $k\equiv 4s+3$ mod $4t$.
Let $$x^*_n =\sum_{4(n-1)t+4s+4\le k\le 4nt+4s+1}P^*_kx^*$$ for $n\ge 0$.
Then for each $n\ge 0$ there exists $x_n\in E(4(n-1)t+4s+4,\allowbreak
4nt+4s+1)$ by
(g) so that $\|x_n\|=1$ and $x^*_n(x_n)\ge (1-\delta)\|x^*_n\|$.  We also
have $x^*_n(x_m)=0$ if $n\neq m$ by condition~(f).

Hence if $(\alpha_n)_{n\ge 0}$ is any finitely nonzero sequence with
$\|(\alpha_n)\|_p=1$, then
$$ \|T_s^*x^*\| \ge 
\left|\left\langle x^{*},\sum \alpha_{n}x_{n} \right\rangle\right| \ge
(1-\delta)\sum \alpha_n\|x^*_n\|.$$
We draw the conclusion that 
$$\biggl(\sum \|x^*_n\|^q\biggr)^{1/q}\le
\frac{1}{1-\delta}\|T_{s}^{*}x^{*}\| \le
\frac{1+\delta}{1-\delta}\|x^*\|.$$

By the remarks above it follows that $\sum_{k=1}^n x_k^*$ is bounded in
norm by $(1+\delta)^2(1-\delta)^{-1}\|x^*\|$.  Clearly $\sum_n
x_n^*(u_k)$ converges for each $k$ and so $\sum x_n^*$ converges weak$^*$
(and actually also in norm, but we do not need this).  Since
$P_i^{*}(x^*-\sum x_n^*)=0$ for all $i$ it follows that $x^*(u_k)=\sum
x_n^*(u_k)$ for all $k$.  Hence $x^*=\sum x_n^*$ weak$^*$.

Now
\bea
\|x^*\| &\le& \sup_n\left\|\sum_{k=1}^nx_k^*\right\|\\
&\le& (1+\delta)\biggl(\sum_{k=1}^{\infty}\|x^*_k\|^q\biggr)^{1/q}\\
&\le& \frac{1+\delta}{1-\delta}\|T^*_sx^*\|
\eea
provided, of course, that $R_s^{*}x^*=0$.

Now let us relax this assumption.  In general the norm of $x^*$
restricted to the range of $R_s$ is at most
$(1-\delta)^{-1}\|R_s^*x^*\|$.
Hence by the Hahn-Banach
theorem there exists $y^*\in X^*$ with $\|y^*\|\le
(1-\delta)^{-1}\|R_s^*x^*\| \le 4\|R_s^*x^*\|$ and
so that $R_s^{*}y^*=R_s^{*}x^*$.  But by the above
$$ \|y^*-x^*\| \le \frac{1+\delta}{1-\delta}\|T_s^*(y^*-x^*)\|.$$
Hence
\bea
\|T_s^*x^*\| &\ge& \|T_s^*(y^*-x^*)\|-\|T_s^*y^*\|  \\
             &\ge& \frac{1-\delta}{1+\delta}\|y^*-x^*\| -
(1+\delta)\|y^*\|\\
   &\ge& \frac{1-\delta}{1+\delta}\|x^*\| - 4\|y^*\|\\
&\ge& \frac{1-\delta}{1+\delta}\|x^*\| -
16\|R_s^*x^*\|.
\eea
{}From this we deduce using the triangle inequality in $\llq^{t}$
that
$$ \biggl(\frac1t\sum_{s=0}^{t-1}\|T_s^*x^*\|^q\biggr)^{1/q} \ge
\frac{1-\delta}{1+\delta}\|x^*\| -
16\biggl(\frac1t\sum_{s=0}^{t-1}\|R_s^*x^*\|^q\biggr)^{1/q}.$$
Hence
$$ \|T^*x^*\| \ge \frac{1-\delta}{1+\delta}\|x^*\|
-\frac{16}{t^{1/q}}\|R^*x^*\|,$$
and this gives
$$ \|T^*x^*\| \ge (1-{\textstyle \frac14}\epsilon)\|x^*\|.$$
Since $\|T\|\le 1+\frac14\epsilon$ this implies that
$d(X,Y/\ker T)< 1+\epsilon$.
Since $Y$ is an $\ell_p$-sum of finite-dimensional spaces we are done.
\eop
%
%
%	Theo F3
%
%
\begin{theo}\label{F3} 
Suppose $1<p<\infty$ and $X$ is a separable Banach
space not containing a copy of $\ell_1$.  
Then $X$ has property $(m_p)$ if and only if, given any
$\epsilon>0$, there is a subspace
$X_0$ of $C_p$ so that  $d(X,X_0)<1+\epsilon$.
\end{theo}
\Proof
In this case $X$, being a quotient of $C_p$ by the preceding
result, is reflexive and so Theorem~\ref{F2} can be applied
to $X^*$ to give the result.
\eop
\bigskip

The following corollary is  a more precise statement of a classical
result of Johnson and Zippin \cite{JohZip2} who have proved the
corresponding isomorphic results.
%
%
%	Cor F4
%
%
\begin{cor}\label{F4}
Let $1<p<\iy$ and $\eps>0$.
\begin{statements}
\item
If $X$ is a quotient of $C_p$,
then there is a subspace $X_0$ of $C_p$ with $d(X,X_0)<1+\epsilon$.
\item
If $X$ is a subspace of $C_p$,
then there is a quotient  $X_0$ of $C_p$ with $d(X,X_0)<1+\epsilon$.
\end{statements}
\end{cor}
\Proof
(a)
We have that $X^*$ has property $(m_q)$, since it is a subspace of
$C_q$ where as usual $1/p+1/q=1$, and hence $X^*$ has		$(m_p^*)$ 
by Theorem~\ref{DA}. It follows that $X$ has
$(m_p)$, and it remains to  apply Theorem~\ref{F3}.

The proof of (b) is similar.
\eop
%
%
%		Theo F5
%
%
\begin{theo}\label{F5}
Let $X$ be a separable Banach space not containing a copy
of $\ell_1$.  Then $X$ has
property $(m_{\infty})$ if and only if, for any $\epsilon>0$, there is a
subspace $X_0$ of $c_0$ with $d(X,X_0)<1+\epsilon$.
\end{theo}
\Proof
If $X$ has $(m_{\infty})$, then Theorem~\ref{F2} shows that there is
a quotient of $C_{\infty}$, $X_1$ say, with $d(X,X_1)<1+\epsilon$.
However, this implies that there is a quotient of a subspace of $c_0$
and hence a subspace of a quotient of $c_0$,
$X_2$ say, with $d(X,X_2)<1+\epsilon$.  The result now follows from the
Alspach's theorem \cite{Als} that quotients of $c_0$ are almost isometric to
subspaces of $c_0$.
\eop
\bigskip

We finally combine results from \cite{Kal-M}, \cite{OjaDi} and 
Theorems \ref{Dtheo}  and \ref{F3} (resp.~\ref{F5})
to obtain the following corollary. If $X$ satisfies
condition~(i) below, it was said to be an $(M_p)$-space in 
\cite{OjaDi}; see also Chapter~VI.5
in \cite{HWW} for information on $(M_p)$-spaces.
%%
%
%		Cor F6
%
%
\begin{cor}\label{F6}
Let $1<p\le\iy$ and let $X$ be a separable Banach space. Then the 
following conditions are equivalent.
\begin{aequivalenz}
\item
$\K(X\oplus_p X)$ is an $M$-ideal in $\LL(X\oplus_p X)$.
\item
The space $X$ has the metric compact approximation property, 
enjoys property $(m_p)$ and fails to contain a copy of $\ell_1$.
\item
The space $X$ has the metric compact approximation property
and is almost isometric to a subspace of an $\ell_p$-sum of
finite-dimensional spaces when $p<\iy$, respectively, to a subspace
of $c_0$ when $p=\iy$. 
\end{aequivalenz}
\end{cor}
\noindent
{\em Remarks.} (1) The equivalence 
(i) $\Longleftrightarrow$ (iii) has been conjectured by $M$-idealists
for some time, cf.\ \cite[p.~336]{HWW}.
A special case of the equivalence 
(i) $\Longleftrightarrow$ (iii) for $p=\iy$
was proved by W.~Werner \cite{WW4}.
Actually, it now turns out that the condition considered by him
is equivalent to the $(M_\iy)$-condition.

(2) One can prove the implication (ii) $\Rightarrow$	(iii)
employing techniques as in the proof of Theorem~\ref{EA}. 
For simplicity suppose $p<\iy$.
In this case one concludes first that there is a sequence $(S_n)$
in $\K(X)$ such that
$$
(1-\eps) \|x\| \le \biggl(\sum_{n=0}^\iy \|S_nx\|^p\biggr)^{1/p}
\le (1+\eps)\|x\|
$$
for all $x\in X$.	Then one embeds $X$ into a space $Y$ with the 
approximation property, say via an embedding $j$. Then $jS_n$
can be approximated up to 		$\eps/2^n$ by finite-rank operators
$F_n$. Thus $X$ embeds into $\ell_p(\ran F_n)$.

(3) For examples illustrating Corollary~\ref{F6} we refer to
Corollaries \ref{G8} and \ref{G9} below.
%\bigskip
%
%
% Section 5
%
%
\mysec{5}{Property $(M)$ for subspaces of $\Lp$}%
The subject-matter of the present section is to refine the results
of the preceding one for subspaces of $\Lp$. 
Here $\Lp$ stands for a separable nonatomic $\Lp(\mu)$-space; it
is known that such a space is isometric to $\Lp[0,1]$ 
\cite[p.~321]{Royden}.
We are first going to prove a  technical
result, Theorem~\ref{G3}, that will be applied to subspaces of
$\Lp$ in this section and to subspaces of the Schatten classes
in the following one. In order to formulate it we need the notions
of a finite-dimensional decomposition (FDD) and of cotype for which we 
refer to \cite[p.~48]{LiTz1} and \cite[p.~72ff.]{LiTz2}, 
respectively. We also found it convenient
to introduce the following technical definition.
%
%
%		Def G2
%
%
\begin{df}\label{G2} \em
Let $Y$ be a separable Banach space and let $\tau$ be a 
topology on $Y$ making it a topological vector space.  
If $1<p<\infty$ we will say that $Y$ has {\em property\/}
$(m_p(\tau))$ if whenever $\|u\|=1$ and $(y_n)$ is a normalized sequence
with $\lim y_n=0$ for the topology $\tau$, then 
$$
\lim_{n\to\iy}\|\alpha u+\beta y_n\|^p=|\alpha|^p+|\beta|^p
$$ 
for all scalars $\alpha$ and $\beta$.
\end{df}

We remark that $Y$ has $(m_{p}(\tau))$ if and only if
$$
\limsup\| u+ y_n\|=\|(\|u\|, \limsup \|y_{n}\|)\|_{p}
$$ 
whenever $\tau$-$\lim y_{n}=0$. In the proof of Proposition~\ref{B4}
we took advantage of the fact $\Lp$ has $(m_{p}(\tau))$ for the
topology of convergence in measure.
%
%
%	Theo G3
%
%
\begin{theo}\label{G3}   
Suppose $Y$ is a separable
Banach space with
nontrivial cotype.  Suppose $Y$ has an unconditional (FDD) 
and denote the partial sum
operators by $S_n$.  Suppose further that
$\tau$ is  a vector topology on $Y$ weaker than the norm
topology and such that $Y$ has property $(m_p(\tau))$ for some $p>1$.

Now let $X$ be a closed infinite-dimensional subspace of\/ $Y$ so that
$S_nx\to x$ for the topology $\tau$ uniformly for $x\in B_X$.
 Then, given $\epsilon>0$ there are a sequence of
finite-dimensional subspaces $(E_n)$ of $Y$ and  a subspace $X_0$ of
$\ell_p(E_n)$ so that $d(X_0,X)\le 1+\epsilon$.
\end{theo}
\Proof 
The argument is a variant of Theorem~\ref{EA}.  
We suppose $Y$ has cotype~$q$ 
where $q\ge p$, with constant $C_q$, and we denote
the unconditional constant of the (FDD) by $C_u$.
  We select
$0<\delta<\frac18\epsilon$ and choose $\eta_n>0$ so that
$\sum\eta_n<\delta$, and finally  we 
pick an integer $t$ so that $C_uC_qt^{-1/q}\le
\frac18\epsilon$.

For convenience we index $(S_n)$ from zero with $S_0=0$.
We will construct a subsequence $T_n=S_{r_n}$ (for $n\ge 0$) of $(S_n)$
with the property that if $A_n=\{\sum_{j=1}^n\lambda_j(T_jx-T_{j-1}x)\dopu
|\lambda_j|\le 1$, $ x\in B_X\}$ and 
$D_n=\{S_mx-S_lx \dopu x\in B_X$, 
$m>l\ge r_n\}$, then for any $v\in D_n$, $u\in A_{n-1}$ we have:
$$ (\|u\|^p+\|v\|^p)^{1/p}-\eta_n \le \|u+v\| \le
(\|u\|^p+\|v\|^p)^{1/p}+\eta_n.$$

To start the induction set $r_0=0$.
Now suppose $r_0,\ldots,r_{n-1}$ have been chosen.  If no successful
choice of $r_n>r_{n-1}$ can be made, then we can find sequences
$r_{n-1}<a_1<b_1<a_2<\cdots$ and $x_k\in B_X$ so that if
$v_k=S_{b_{k}}x_k-S_{a_k}x_k$, then there exist $u_k\in A_{n-1}$ with
$$ \bigl|\,\|u_k + v_k\| -(\|u_k\|^p+\|v_k\|^p)^{1/p}\,\bigr|>\eta_n.$$
Passing to a subsequence we can suppose that $(u_k)$ converges in norm to
some
$u\in A_{n-1}$.  However $(v_k)$ then converges to $0$ for $\tau$
by assumption on $X$,
and we obtain from property~$(m_p(\tau))$
 $$ \lim_{k\to \infty}(\|u+v_k\|-(\|u\|^p+\|v_k\|^p)^{1/p})=0$$
which leads to a contradiction.
Thus the induction can be carried out.

Now it follows by induction that, if $x\in B_X$ and $(a_k)_{k\ge 0}$,
$(b_k)_{k\ge 0}$ are two sequences satisfying $0\le
a_0<b_0<a_1<b_1<\cdots$, then, letting  $v_k=T_{b_k}x-T_{a_k}x$,
we have
$$ 
\biggl(\sum_{j=1}^n \|v_j\|^p\biggr)^{1/p} -\sum_{j=1}^n\eta_{a_j} \le
\left\|\sum_{j=1}^nv_j \right\| \le
\biggl(\sum_{j=1}^n \|v_j\|^p\biggr)^{1/p} +\sum_{j=1}^n\eta_{a_j}.
$$
Since $\sum v_j$ converges as the (FDD) is unconditional this means that
$$ 
\biggl(\sum_{j=1}^{\infty} \|v_j\|^p\biggr)^{1/p} -\delta \le
\left\|\sum_{j=1}^{\infty}v_j \right\| \le
\biggl(\sum_{j=1}^{\infty} \|v_j\|^p \biggr)^{1/p} +\delta.
$$

Now let $V_k=T_{k+1}-T_k$ for $k\ge 0$, and for $0\le s\le t-1$ define
$A_sy= \sum_{j\equiv s}V_jy$ (congruence mod $t$).  We also let
$$W_{k,s}=
\sum_{(k-1)t+s<j<kt+s}V_j$$ for $k\ge 0$ and $0\le s\le t-1$.

Then for any $|\lambda_s|=1$ and any $x\in X$ with $\|x\|=1$
we have
$$ \left\|\sum_{s=0}^{t-1}\lambda_sA_sx\right\| \le C_u$$
and hence by the cotype $q$ condition
$$ \biggl(\sum_{s=0}^{t-1}\|A_sx\|^q\biggr)^{1/q} \le C_uC_q.$$
{}From this we deduce since  $p\le q$,
$$ \biggl(\sum_{s=0}^{t-1}\|A_sx\|^p\biggr)^{1/p} \le C_uC_q t^{1/p-1/q}.$$

Now
$$ 
\biggl(\sum_{k=0}^{\infty}\|W_{k,s}x\|^p\biggr)^{1/p}-\delta 
\le \|x-A_sx\| \le
\biggl(\sum_{k=0}^{\infty}\|W_{k,s}x\|^p\biggr)^{1/p}+\delta.
$$
Hence
$$ 
1-\|A_sx\| -\delta \le 
\biggl(\sum_{k=0}^{\infty}\|W_{k,s}x\|^p\biggr)^{1/p}\le
1+\|A_sx\|+\delta.
$$
It follows that
$$ 
t^{1/p}(1-\delta -C_uC_qt^{-1/q})\le
\biggl(\sum_{k=0}^{\infty}\sum_{s=0}^{t-1}\|W_{k,s}x\|^p\biggr)^{1/p} \le
 t^{1/p}(1+\delta +C_uC_qt^{-1/q}).
$$
By choice of $\delta$ and $t$ this implies that 
$$
1-{\textstyle \frac14\epsilon} \le 
\biggl(\frac1t\sum_{k=0}^{\infty}\sum_{s=0}^{t-1}\|W_{k,s}x\|^p\biggr)^{1/p} 
\le 1+{\textstyle \frac14\epsilon}.
$$

Hence $X$ has Banach-Mazur distance less than $1+\epsilon$ from a closed
subspace of $\ell_p(E_{k,s})$ where $E_{k,s}=\ran W_{k,s}$. 
\eop
\bigskip

It remains to give concrete examples.
We consider $L_p=L_p(\mu)$ where $\mu$ is a
probability measure on a compact metric space $\Omega$
and  $1<p<\infty$.
%
%
%	Lemma G4
%
%
\begin{lemma}\label{G4} 
Suppose $1<p<\infty$ and $p\neq 2$.  Suppose $f\in L_p$
and
$(g_n)$ is a bounded sequence in
$L_p$.   Then the following statements are
equivalent:
\begin{aequivalenz}
\item
For any scalars $\alpha$ and $\beta$ we have
$$ \lim_{n\to\infty}\bigl(\|\alpha f+\beta g_n\|^p-
   |\alpha|^p\|f\|^p-|\beta|^p\|g_n\|^p\bigr)=0.$$
\item
\hfil 
$\displaystyle \lim_{n\to\infty}\int_{|f|>0}|g_n|\,d\mu=0$.
\hfill
\end{aequivalenz}
\end{lemma}
\Proof
(i) $\Rightarrow$ (ii):  
It suffices to consider the case
$\|f\|_p=\|g_n\|_p=1$.  Then by a result of Dor \cite{Dor} 
(cf.\ the discussion in
Alspach \cite{Als2}), there exist for each $n$ disjoint Borel sets 
$A_n$, $B_n$ so
that $\|f\chi_{A_n}\|_p\to 1$ and $\|g\chi_{B_n}\|_p\to 1$.  Hence
$\lim\|g\chi_{A_n}\|_p=0$.  Now let $C_n=\{\omega\dopu |f(\omega)|>0$,
$\omega\notin A_n\}$.  Then $\lim\mu(C_n)=0$.  Hence $(g_n)$
tends to~0 in measure
on $\{|f|>0\}$ and is equi-integrable, so that (ii) follows.

(ii) $\Rightarrow$ (i):  This is standard and we omit it.
\eop
%
%
%		Theo G5
%
%
\begin{theo}\label{G5}  
Suppose $1<p<\infty$, $p\neq2$,  and let $X$ be a closed
infinite-dimen\-sional subspace of
$L_p$.  Then the following conditions on $X$ are equivalent:
\begin{aequivalenz}
\item
$B_X$ is compact in $L_1$.
\item
$X$ has property $(m_p)$.
\item
For any $\epsilon>0$, there is a subspace $X_0$ of $\ell_p$ so that
$d(X_0,X)<1+\epsilon$.
\end{aequivalenz}
\end{theo}
\Proof
(iii) $\Rightarrow$ (ii):  This is immediate.

(ii) $\Rightarrow$ (i): 
We  start by producing a function $f$ in $X$
with maximal support.  Let $(f_n)$ be a dense sequence in the unit ball
of $X$.  Then $\sum_{n=1}^{\infty}t_n2^{-n}f_n$ converges in $L_p$ and
absolutely a.e.\ for every $(t_n)\in [0,1]^{\N}$.  By a standard
Fubini argument there is a choice of $t=(t_n)$ so that if $f=\sum
t_n2^{-n}f_n$ then $f$ is nonzero a.e.\ on the set 
$S=\bigcup \{|f_n|>0\}$. (In fact, for almost every $\omega\in S$ the
set of those $t$ for which $\sum t_{n}2^{-n}f_{n}(\omega)=0$ has
measure~0, and for $A=\{(t,\omega)\dopu \omega\in S$, $\sum
t_{n}2^{-n}f_{n}(\omega)\neq0\}$, 
$\int \!\!dt \int \!\!d\mu\, \chi_{A} 
= \int \!\!d\mu \int \!\!dt \,\chi_{A} =
\mu(S)$.)
It is then easy to see that  any $g\in X$ vanishes
a.e.\ on the set $\{f=0\}$.

Now suppose $(g_n)$ is a weakly null sequence in $X$.  Then by $(m_p)$
and the preceding lemma, $\lim_{n\to\infty}\int_{|f|>0} |g_n|\,d\mu=0$.
Thus $\lim\|g_n\|_1=0$ and, as $X$ is reflexive, the inclusion 
$X\hookrightarrow L_1$
is compact and $B_X$ is therefore compact in $L_1$.

(i) $\Rightarrow$ (iii):  
We let $\tau$ be the $L_1$-norm topology on $L_p$.
Then $L_p$ has $(m_p(\tau))$ and nontrivial cotype.  
Also $L_p$ has an unconditional basis (the
Haar basis) \cite[Th.~2.c.5]{LiTz2}.  
The Haar basis is also a basis for $L_1$ \cite[Prop.~2.c.1]{LiTz2}.  
Hence if $(S_n)$
denotes the sequence of partial sum operators, 
then $S_nf\to f$ $L_1$-uniformly on
$B_X$.  We are thus in a position to apply Theorem~\ref{G3}, 
and since each $E_n$ produced there is a finite-dimensional subspace of 
$L_p$ we obtain (iii).
\eop
\bigskip

Next
we note that in this context property~$(M)$ necessarily reduces to
property~$(m_r)$ for some $r$.
%
%
%		Cor G1
%
%
\begin{cor}\label{G1}
Suppose $1\le p<\infty$ and let $X$ be a closed
infinite-dimen\-sional subspace of $L_p$ with property~$(M)$
and, if $p=1$, failing to contain a copy of $\ell_1$.
Then
there exists $1<r<\infty$ so that $X$ has $(m_r)$, and
for any $\epsilon>0$ there exists a
closed subspace $X_0$ of $\ell_r$ with $d(X,X_0)<1+\epsilon$.
If $p\le 2$, then $p\le r\le2$; if $p>2$, then $r=2$ or $r=p$.
\end{cor}
\Proof
Since $X$ is stable \cite{KriMau} it follows that $X$ has property
$(m_r)$ for some $1<r<\infty$. (This is what the proof of 
\cite[Th.~3.10]{Kal-M} actually shows.)  Theorem~\ref{F3} implies
that $X$ contains a copy of $\ell_{r}$. Now
$\Lp$ contains only those $\ell_{r}$ as stated
\cite[p.~132ff.]{LiTz}, hence we get the final assertion.

 Since there is a weakly null sequence $(u_n)$ in $X$
with $\lim\|x+\alpha u_n\|=(\|x\|^r+|\alpha|^r)^{1/r}$ for any $x\in X$
\cite[Prop.~3.9]{Kal-M},
the space $X\oplus_r\R$ is finitely representable in $X$ and thus
embeds isometrically into an ultraproduct of
$X$ and hence also into $L_p$.  
We wish to conclude that $X$ is isometric to a subspace of $L_{r}$.

If $p\le2$, this is a result independently due to Dor \cite{Dor2} and 
Raynaud  \cite[Prop.~3]{Ray2}. If $p>2$ and $r=2$, we observe
that Raynaud's argument still works if $p$ is not an even integer;
and the assertion is trivial if $r=p$. Suppose now that $p=2m$ is an
even integer and $r=2$. We claim that a Banach space $Y$ is isometric
to a Hilbert space provided $Y\oplus_{2}\R$ is isometric to a
subspace of $\Lp$. Indeed, if $Y\subset \Lp$ and $g\in \Lp$ satisfy
$\|y+tg\|_{p}= (\|y\|_{p}^{2} + |t|^{2})^{1/2}$ for all $y\in Y$,
$t\in\R$, then
$$
\|y+tg\|_{p}^{p} = \int |y+tg|^{2m}d\mu = (\|y\|_{p}^{2}+ |t|^{2})^{m}
$$
is  a polynomial in $t$, the coefficient of $t^{2m-2}$ being
$$
m(2m-1)\int |y|^{2}|g|^{2m-2}d\mu = m\|y\|_{p}^{2}.
$$
Hence $Y$ is a Hilbert space.

To prove the corollary it is now left to apply Theorem~\ref{G5}.
\eop
\bigskip

The case $p>2$ is especially interesting.
%
%
%		Cor G6
%
%
\begin{cor}\label{G6}
Suppose $p>2$ and let $X$ be a closed subspace of
$L_p$.  Then the following conditions on $X$ are equivalent:
\begin{aequivalenz}
\item
$X$ contains no subspace isomorphic to $\ell_2$.
\item
$X$ is isomorphic to a subspace of $\ell_p$.
\item
For any $\epsilon>0$, there is a subspace $X_0$ of $\ell_p$ so that
$d(X_0,X)<1+\epsilon$.
\end{aequivalenz}
\end{cor}
\Proof
The only step requiring proof is that (i) implies (iii). By 
Theorem~\ref{G5} it remains to show that $B_X$
is $L_1$-compact.  If not, there is a sequence $(f_n)$ which is weakly
null but does not converge to zero in $L_1$.  By passing to a
subsequence we can assume that $(f_n)$ is an unconditional basic
sequence, cf.\ \cite[Prop.~1.a.11]{LiTz1}, and that $\inf
\|f_{n}\|_{1}>0$.
We thus have for any finite sequence of scalars $a_1,\ldots,a_n$
\bea
\left\|\sum_{k=1}^n a_kf_k\right\|_p &\ge& 
c_{1}\, {\rm ave}\left\|\sum_{k=1}^n {\pm a_k}f_k\right\|_p\\
&\ge&
c_{2} \left\|\biggl(\sum_{k=1}^n|a_k|^2|f_k|^2\biggr)^{1/2}\right\|_p\\
&\ge& 
c_{2} \left\|\biggl(\sum_{k=1}^n|a_k|^2|f_k|^2\biggr)^{1/2}\right\|_2\\
&\ge& 
c_{2} \inf \|f_{n}\|_{2} \biggl(\sum_{k=1}^n|a_k|^2\biggr)^{1/2}\\
&\ge&
c_{3}\biggl(\sum_{k=1}^n|a_k|^2\biggr)^{1/2}
\eea
for suitable constants $c_{i}>0$.  Since $L_p$ has type 2, it follows that
conversely
\bea
\left\|\sum_{k=1}^n a_kf_k\right\|_p &\le& 
c_{4}\, {\rm ave}\left\|\sum_{k=1}^n {\pm a_k}f_k\right\|_p\\
&\le&
c_{5}\biggl(\sum_{k=1}^n|a_k|^2\biggr)^{1/2}.
\eea
Thus
$\linq\{f_1,f_2,\ldots\}$ is isomorphic to a Hilbert space.
\eop
\bigskip

\noindent
{\em Remarks.}
(1)
The equivalence of (i) and (ii) above, however with a large
isomorphism constant, is due to Johnson and
Odell \cite{JohOde}; see also \cite{John3}.
%as pointed out to us by W.~B.~Johnson, a
%somewhat better constant is obtained in \loch.

(2)
In the case $1<p<2$ the equivalence of (ii) and (iii) breaks down.  In fact
the space $\ell_p\oplus_p\ell_2$ is isometric to a subspace of $L_p$ and
has a subspace isometric to the space $\ell_p$ equipped with the norm
$\|\xi\|=(\|\xi\|_p^p +\|\xi\|_2^p)^{1/p}$.  
A direct calculation shows that this norm does not
 have property $(m_p)$ (not even property $(M)$); so $\lp$
with this norm cannot embed almost isometrically into $\lp$
with its classical norm. 
(A  more roundabout way to arrive at the same conclusion is to
observe that the above norm is 1-symmetric and to
invoke a result due to  Hennefeld \cite{Hen2}  
(see \cite[Prop.~VI.4.24]{HWW} for a proof based on property $(M)$)
according to which the spaces $(\lp,\|~\|_{p})$
are the only 1-symmetric Banach spaces for which $\K(X)$
is an $M$-ideal in $\LL(X)$	and thus the only ones with 
property $(M)$.)
%
%
% Cor G7
%
%
\begin{cor}\label{G7}
Let $2<p<\iy$ and suppose $X$ is a closed infinite-dimen\-sional
subspace of $\Lp$ with property~$(M)$. Then either $X$ embeds almost
isometrically into $\lp$, or $X$ is isometric to a Hilbert space.
\end{cor}
\Proof
This follows from Corollaries \ref{G1} and \ref{G6}. (Note that if $X$
embeds almost isometrically into $\ell_{2}$, then the parallelogram
identity holds, and $X$ is isometric to a Hilbert space.)
\eop
\bigskip

As a specific example of a subspace of $\Lp$ we now discuss
the Bergman space~$B_p$.
This is the space of all analytic functions on the open unit disk $\D$
in the complex plane for which $\int_{\D}|f(x+iy)|^{p}\,dxdy$ is
finite. Thus $B_{p}$ is a linear subspace of $\Lp=\Lp(\D,dx\,dy)$,
and it is closed. It is well known that $B_{p}$ has the metric
approximation property (use Fej\'{e}r kernels, as in
\cite[Lemma~3.4]{Dirk5}); in fact, it is isomorphic to $\lp$
\cite{LinPel2}. Let us point out that $B_{p}$ has $(m_{p})$. If
$f_{n}\to 0$ weakly, then clearly $f_{n}\to 0$ pointwise. Since
bounded sets in $B_{p}$ are equicontinuous on compact subsets of
$\D$, we conclude that $f_{n}\to 0$ uniformly on compact subsets of
$\D$. From this it follows easily that $B_{p}$ has $(m_{p})$.

Hence Corollary~\ref{F6} and Theorem~\ref{G5} imply:
%
%
% Cor G8
%
%
\begin{cor}\label{G8}
The Bergman space $B_{p}$ embeds almost isometrically into $\lp$,
and $\K(B_{p})$ is an $M$-ideal in $\LL(B_{p})$ if $1<p<\iy$.
\end{cor}

We also wish to consider the corresponding space of analytic
functions when $p=\iy$.  This is the so-called `little' Bloch space
$\beta_{0}$ consisting of those analytic functions $f$ on $\D$ for
which $\|f\|_{\beta}:= \sup (1-|z|^{2})|f(z)|<\iy$ and
$\lim_{|z|\to1}(1-|z|^{2})f(z)=0$. (The analytic functions satisfying
$\|f\|_{\beta}<\iy$ form the Bloch space $\beta$. We remark for
completeness that in
complex analysis the space of functions whose derivative belongs to
$\beta$ is called the Bloch space.) It is known that
$\beta$ is the bidual of $\beta_{0}$ (e.g., \cite{BieSum}), $\beta_{0}$
is an $M$-ideal in $\beta$ \cite{Dirk5}, $\beta_{0}$ has $(m_{\iy})$
(the above argument works here as well, see also \cite{AETD}), and
$\beta_{0}$ has the metric approximation property (e.g.,
\cite[Lemma~3.4]{Dirk5}). Hence we get from Corollary~\ref{F6}:
%
%
% Cor G9
%
%
\begin{cor}\label{G9}
The little Bloch space $\beta_{0}$ embeds almost isometrically into
$c_{0}$,
and $\K(\beta_{0})$ is an $M$-ideal in $\LL(\beta_{0})$.
In fact, for every Banach space $X$, $\K(X,\beta_{0})$ 
is an $M$-ideal in $\LL(X,\beta)$.
\end{cor}

This result (the last mentioned assertion follows from \cite{Dirk2})
was conjectured in \cite{AETD}. Note that $\beta_{0}$ is not
isometric to a subspace of $c_{0}$, because its unit ball has extreme
points \cite{CimWo1}, but the only subspaces of $c_{0}$ whose unit
balls have extreme points are finite-dimen\-sional.

Actually, to obtain the above  corollaries it is not necessary to
apply the machinery developed in this paper; all one has to know is
that $B_{p}$ embeds almost isometrically into $\lp$ and $\beta_{0}$
into $c_{0}$. This can be proved directly as follows; for
simplicity let us give the details for $\beta_{0}$. Let $\eps>0$ and
$r_{n}= 1-2^{-n}$ for $n\ge0$. Then for each $n$ there are constants
$K_{n}$ and $L_{n}$ so that if $\|f\|_{\beta}=1$, then 
$|f(z)|\le K_{n}$ and $|f'(z)|\le L_{n}$ for
$|z|\le r_{n}$. Now for each $n$ choose a finite subset $F_{n}$ of
the annulus $A_{n}=\{z\dopu r_{n-1}\le |z| \le r_{n}\}$ so that if
$z\in A_{n}$, then there exists $z'\in F_{n}$ with $|z-z'|\le
\delta_{n}$  where $(L_{n}+2K_{n})\delta_{n}\le\eps$.

Suppose $\|f\|_{\beta}=1$. If $z\in A_{n}$, pick $z'\in F_{n}$ as
above. Hence $|f(z)|\le |f(z')|+L_{n}\delta_{n}$ and so
$$
(1-|z|^{2})|f(z)| \le 
(|z'|^{2}-|z|^{2})|f(z)| + (1-|z'|^{2})(|f(z')|+L_{n}\delta_{n}).
$$
Hence
$$
\sup_{z\in A_{n}} (1-|z|^{2})|f(z)| \le
(L_{n}+2K_{n})\delta_{n} + \max_{z\in F_{n}} (1-|z|^{2})|f(z)| .
$$
Thus if $F=\bigcup F_{n}$,
$$
1 \le \sup_{z\in F} (1-|z|^{2})|f(z)| + \eps.
$$
If we arrange the set $F$  into a sequence $(z_{n})$, then the
operator $T$ defined by 
$Tf= \bigl((1-|z_{n}|^{2})f(z_{n})\bigr)$ maps $\beta_{0}$ into
a subspace of $c_{0}$ and $(1-\eps)\|f\|_{\beta}\le \|Tf\| \le
\|f\|_{\beta}$.

Let us remark that our arguments are not restricted to the weight
${1-|z|^{2}}$ appearing in the Bloch norm; as a matter of fact one may
consider any radial weight of the form $v(|z|)$ where $v\dopu
[0,1]\to [0,1]$ is a continuous decreasing function with $v(r)=0$ if
and only if $r=1$. For all these weighted spaces of analytic
functions Corollary~\ref{G9} is valid. Lusky (\cite{Lusk2},
\cite{Lusk3}) has recently investigated those spaces and shown that
they always embed into $c_{0}$ isomorphically, and he has determined
for which weights one actually gets spaces isomorphic to $c_{0}$.
%\bigskip
%
%
% Section 6
%
%
\mysec{6}{Property $(M)$ for subspaces of the Schatten classes}%
At this point we also consider the Schatten classes 
$c_p$ of operators on a separable Hilbert space $\cal H$
and prove  similar
results like for $\Lp$ by essentially the same technique.  
The reader is cautioned that $x^*$ will
now denote the Hilbert space adjoint of an operator $x\in c_p$.
We will need the following analogue of Lemma~\ref{G4} above. We use
$\|~\|_{\infty}$ to denote the operator norm.
%
%
%		Lemma H2
%
%
\begin{lemma}\label{H2}
Suppose $1<p<\infty$ and $p\neq 2$.  Suppose $x\in c_p$
and $(y_n)$ is a bounded sequence in $c_p$.  The following conditions are
equivalent:
\begin{aequivalenz}
\item
For any scalars $\alpha$ and $\beta$ we have
$$ \lim_{n\to\infty}\bigl(\|\alpha x+\beta
y_n\|_p^p-|\alpha|^p\|x\|_p^p-|\beta|^p\|y_n\|_p^p\bigr)=0.$$
\item
\hfil
$ \displaystyle \lim_{n\to\infty}\|y_n^*x\|_{\infty}=\lim_{n\to\infty}
\|y_nx^*\|_{\infty} =0$.
\hfill
\end{aequivalenz}
 \end{lemma}
\Proof
(ii) $\Rightarrow$ (i): 
Given $\epsilon>0$ we can choose
finite-rank
orthogonal projections $\pi$, $\pi'$ so that $\|\pi'x\pi-x\|_p<\epsilon$
and $\ran (\pi) \subset \ran (x^*)$, $\ran (\pi')\subset \ran (x)$.
Now $y_n(\xi)\to 0$ for all $\xi\in \ran (x^*)$ and so
$\lim_{n\to\infty}\|y_n\pi\|_{\infty}=0$.  Since $\pi$ has
finite rank,
$\lim_{n\to\infty}\|y_n\pi\|_p=0$.   Similarly
$\lim_{n\to\infty}\|\pi'y_n\|_p=0$.

Now it follows that for any $\alpha$ and $\beta$,
$$
\lim_{n\to\infty}
\bigl(\|\alpha \pi' x\pi+\beta y_n\|_p -\|\alpha \pi'x\pi
+\beta(\Id-\pi')y_n(\Id-\pi)\|_p\bigr) = 0.$$
Since  $$\|\alpha \pi'x\pi +\beta
(\Id-\pi')y_n(\Id-\pi)\|_p^p=|\alpha|^p\|\pi'x\pi\|_p^p
+|\beta|^p\|(\Id-\pi')y_n(\Id-\pi)\|_p^p$$
we quickly obtain that
$$\limsup_{n\to\infty}\bigl|\,\|\alpha x +\beta
y_n\|_p-(|\alpha|^p\|x\|^p+|\beta|^p\|y_n\|_p^p)^{1/p}\,\bigr| \le
|\alpha|\epsilon.$$
Since $\epsilon>0$ is arbitrary we have the result.

(i) $\Rightarrow$ (ii): Let us start by calling a sequence $(y_n)$ in $c_p$
{\em tight\/} if the induced sequences of singular values
$\sigma_n=(s_k(y_n))_{k=1}^{\infty}$ form a relatively compact subset of
$\ell_p$.  This is easily seen to be equivalent to the requirement that
given any $\epsilon>0$ there exists a natural number $t=t(\epsilon)$ so
that for each $n\in\N$ there is an orthogonal projection $\pi_n$ of rank at
most $t$ so that $\|y_n-y_n\pi_n\|_p<\epsilon$.

Let us first assume that $(y_n)$ is a tight sequence.  Let $\CalU$
be any non-principal ultrafilter on the natural numbers and consider the
ultraproduct $\CalH_{\CalU}$ defined by $\ell_{\infty}(\CalH)$
quotiented by the subspace of all sequences $(\xi_n)$ so that
$\lim_{n\in\CalU}\|\xi_n\|=0$.

We will consider the maps $\theta\in\LL(\CalH_{\CalU})$ defined by
$\theta\bigl((\xi_n)\bigr)=\bigl(x(\xi_n)\bigr)$ and $\phi$ defined by
$\phi\bigl((\xi_n)\bigr)=\bigl(y_n(\xi_n)\bigr)$.  
Suppose $\alpha$ and $\beta$ are any scalars.
It is clear that the singular values of $\alpha \theta+\beta \phi$ are
given by $$s_k(\alpha\theta+\beta\phi)=\lim_{n\in\CalU}s_k(\alpha x +\beta
y_n).$$

We will argue that $(\alpha x+\beta y_n)$ is a tight sequence.  In fact,
given $\epsilon>0$ we can find $t\in\N$ so that for each $n$ there is an
orthogonal projection $\pi_n$ of rank at most $t$ so that
$\|\beta(y_n-y_n\pi_n)\|_p<\epsilon/2$.  We also fix a  
finite-rank projection $\pi_0$ of rank $t_0$ so that
$\|\alpha(x-x\pi_0)\|_p<\epsilon/2$.  Let $\pi'_n$ be the orthogonal
projection onto $\ran (\pi_0)+\ran (\pi_n)$.  Then $\pi'_n$ has rank at
most $t+t_0$ and $\|(\alpha x +\beta y_n)(\Id-\pi'_n)\|_p<\epsilon$.

It thus follows that $\alpha\theta+\beta\phi\in c_p(\CalH_{\CalU})$ and
that $$\|\alpha\theta+\beta\phi\|_p =\lim_{n\in\CalU}\|\alpha x+\beta
y_n\|_p.$$
Hence we have that
$$\|\alpha\theta+\beta\phi\|_p^p=|\alpha|^p\|\theta\|_p^p
+|\beta|^p\|\phi\|_p^p.$$

A classical result of McCarthy \cite{McCar}
gives that in the case of
equality in the generalized Clarkson inequalities in $c_p$ we have
$|\theta||\phi|=0$ or $\phi\theta^*=\phi^*\theta=0$.
Translating this back we obtain that
$$
\lim_{n\to\infty}\|y_nx^*\|_{\infty}=\lim_{n\to\infty}\|y_n^*x\|_{
\infty}  =0.
$$
This completes the proof when $(y_n)$ is tight.

For the general case we may reduce the problem by passing to a
subsequence so that $\lim_{n\to\iy} s_k(y_n)=s_k$ exists for all $k$.  
Then if $y_n$ is
written in its Schmidt series
$y_n=\sum_{k=1}^{\infty}s_k(y_n)\,\xi_k\otimes \eta_k$ where $(\xi_k)$ and
$(\eta_k)$ are orthonormal sequences, we will put
$y'_n=\sum_{k=1}^{\infty}s_k\,\xi_k\otimes \eta_k$.  Then $(y_n')$ is tight
and $\lim\|y_n-y'_n\|_{\infty}=0$.  To see the latter statement notice that
$s_k(y_n)\le Mk^{-1/p}$ where $M=\sup\|y_n\|_p$ and so $|s_k(y_n)-s_k|\le
2Mk^{-1/p}$ but $\lim_{n\to\infty}s_k(y_n)-s_k=0$ for each $k$.

Now suppose $\alpha$ and $\beta$ are given.  Then $(\alpha x+\beta y'_n)$ is
tight as before and hence we can find isometries $u_n$ and $v_n$ so that the
set $(u_n(\alpha x+\beta y'_n)v_n)$ is relatively compact.
However $\|u_n(y_n-y_n')\|_{\infty}\to 0$ and by (ii) $\Rightarrow$ (i) and
a standard compactness argument we get
$$\|u_n(\alpha x +\beta y_n)v_n\|_p^p
= \|u_n(\alpha x +\beta y'_n)v_n\|_p^p +|\beta|^p\|u_n(y_n-y'_n)v_n\|_p^p
+\epsilon_n$$
where $\lim\epsilon_n=0$.  Thus
$$ \|\alpha x+\beta y_n\|_p^p = \|\alpha x+\beta y'_n\|_p^p
+|\beta|^p\|y_n-y'_n\|_p^p+\epsilon_n.$$

We can apply the same equation with $\alpha=0$ to obtain
$$ \lim \bigl( \|y_n\|_p^p-\|y'_n\|_p^p-\|y_n-y'_n\|_p^p \bigr)=0.$$
Putting this together we see that
$$ \lim \bigl(\|\alpha x +\beta y'_n\|^p_p
-|\alpha|^p\|x\|_p^p-|\beta|^p\|y'_n\|_p^p \bigr) =0.$$
Since $(y'_n)$ is tight this means that
$\lim\|y_n'x^*\|_{\infty}=\lim\|(y'_n)^*x\|_{\infty}=0$ and hence
$\lim\|y_nx^*\|_{\infty}=\lim\|y_n^*x\|_{\infty}=0$.
This completes the proof.
\eop
\bigskip

To state our main theorem let $\tau$ be the topology on $c_p$ of
given by the seminorms $x\mapsto \|x(\xi)\|$ and $x\mapsto \|x^*(\xi)\|$ for
$\xi\in\CalH$.  Denote by $c_p^n$ the space of $n\times n$-matrices with
the $c_p$-norm.
%
%
%		Theo H3
%
%
\begin{theo}\label{H3}
Suppose $1<p<\infty$ with $p\neq 2$, and $X$ is a
closed
infinite-dimensional subspace of $c_p$. Then the following conditions are
equivalent:
\begin{aequivalenz}
\item
 $B_X$ is compact for the topology $\tau$.
\item
$X$ has property $(m_p)$.
\item
For any $\epsilon>0$ there is  subspace $X_0$ 
of $\ell_p^{\vphantom{n}}(c^n_p)$ so
that $d(X,X_0)<1+\epsilon$.
\end{aequivalenz}
\end{theo}

\Proof
We note first that $c_p$ has nontrivial cotype (\cite{McCar} and
\cite{TJ}).  Let
$(\xi_n)$ be any fixed orthonormal basis of $\CalH$ and let $\pi_n$ be
the orthogonal projection onto $\lin\{\xi_1,\ldots,\xi_n\}$.  
Then $c_p$ has an unconditional
(FDD) with partial sum operators $S_n(x)=\pi_nx\pi_n$, cf.\
\cite{AraLin}.

We claim that $c_p$ has $(m_p(\tau))$. In fact if $\|x\|_p=1$ and
$\|y_n\|_p=1$ with $y_n\to 0$ for $\tau$, then since $x$ is compact
$\|y_nx^*\|_{\infty} \to 0$ and $\|y_n^*x\|_{\infty}\to 0$.  Hence
$\lim\|\alpha x+\beta y_n\|_p=(|\alpha|^p+|\beta|^p)^{1/p}$
by Lemma~\ref{H2}.

Now we give the proof of the theorem:

(iii) $\Rightarrow$ (ii): Obvious.

(ii) $\Rightarrow$ (i): 
Suppose $(y_n)$ converges to $0$ weakly in $X$.
Then for any $x\in X$ we have by Lemma~\ref{H2}
that $\lim\|y^*_nx\|_{\infty}=0$.  In
particular it follows that $\lim\|y_n^*(\xi)\|=0$ whenever $\xi$ is in the
closed linear span of $\bigcup_{x\in X} \ran (x)$.  From this it follows
easily that $(y_n^*)$ converges to $0$ for the strong operator topology.
Similarly $(y_n)$ converges to 0 for the strong operator topology and so
$y_n\to 0$ for $\tau$ which implies that $B_X$ is $\tau$-compact.

(i) $\Rightarrow$ (iii): 
We need only show that $S_nx\to x$
$\tau$-uniformly for $x\in B_X$.   In fact if $\xi\in \CalH$, then
$\{x(\xi)\dopu x\in B_X\}$ is compact and so
 $$ \limp_{n\to\infty}\sup_{x\in B_X}\|\pi_nx(\xi)-x(\xi)\|=0.$$
On the other hand, if $x\in B_X$, then 
$$\|\pi_nx\pi_n(\xi)-\pi_nx(\xi)\|\le
\|\xi-\pi_n(\xi)\|.$$
 Combining we have that
 $$ \limp_{n\to\infty}\sup_{x\in B_X}\|\pi_nx\pi_n(\xi)-x(\xi)\|=0.$$
This with a similar statement for $x^*$ allows us 
to apply Theorem~\ref{G3} to
obtain the result.
\eop
\bigskip

The following result is analogous to Corollary~\ref{G1}.
%
%
%		Cor H1
%
%
\begin{cor}\label{H1}
Suppose $1< p<\infty$ and let $X$ be a closed
infinite-dimen\-sional subspace of $c_p$ with property~$(M)$.
Then
there exists $r\in\{2,p\}$ so that $X$ has $(m_r)$.
If $X$  does not contain a subspace isomorphic to $\ell_{2}$,
for any $\epsilon>0$ there exists a
closed subspace $X_0$ of $\lpcpn$ with $d(X,X_0)<1+\epsilon$.
\end{cor}
\Proof
We argue as before that $X$ is stable, which was proved by Arazy
\cite{Ara2}. Therefore, for some $r$, $X$ has $(m_{r})$ and contains
a copy of $\ell_{r}$. This implies that $r=2$ or $r=p$ by
\cite{AraLin}. It remains to apply Theorem~\ref{H3}.
\eop
\bigskip

Again, we specialize to $p>2$.
%
%
%	Cor H4
%
%
\begin{cor}\label{H4}
Suppose $2<p<\infty$ and that $X$ is a subspace of
$c_p$. Then the following conditions on $X$ are equivalent:
\begin{aequivalenz}
\item
$X$ contains no  subspace isomorphic to $\ell_2$.
\item
$X$ is isomorphic to a subspace of $\lpcpn$.
\item
For any $\epsilon>0$ there exists a closed subspace $X_0$
of $\lpcpn$ so that $d(X,X_0)<1+\epsilon$.
\end{aequivalenz}
\end{cor}
\Proof
The argument is very similar to that of Corollary~\ref{G6}. We only
have to show that (i) implies (iii), that is, by Theorem~\ref{H3},
that $B_{X}$  is $\tau$-compact. If not, there would be a weakly null
sequence $(x_{n})$ in $B_{X}$ which is an unconditional basic
sequence and satisfies $\inf\|x_{n}(\xi)\|>0$ or 
$\inf\|x^{*}_{n}(\xi)\|>0$ for some $\xi\in\CalH$. (One can extract
an unconditional basic sequence since $c_{p}$, having an
unconditional (FDD), embeds into a space with an unconditional basis
\cite[Th.~1.g.5]{LiTz1}.) Let us assume without loss of generality
that $\inf\|x_{n}(\xi)\|>0$. We then have for suitable positive
constants $c_{i}$ and all finite sequences of scalars
$a_{1},\ldots,a_{n}$
\bea
\left\|\sum_{k=1}^n a_kx_k\right\|_p &\ge& 
c_{1}\, {\rm ave}\left\|\sum_{k=1}^n {\pm a_k}x_k\right\|_p\\
&\ge&
c_{1}\, {\rm ave}\left\|\sum_{k=1}^n {\pm a_k}x_k(\xi)\right\|\\
&\ge&
c_{2} \biggl(\sum_{k=1}^n \| a_k x_k(\xi)\|^{2}\biggr)^{1/2}\\
&\ge&
c_{3}\biggl(\sum_{k=1}^n|a_k|^2\biggr)^{1/2}.
\eea
On the other hand, $c_{p}$ has type~2 (\cite{McCar} and \cite{TJ}),
and we get the reverse estimate
$$
\left\|\sum_{k=1}^n a_kx_k\right\|_p \le
c_{4}\biggl(\sum_{k=1}^n|a_k|^2\biggr)^{1/2}.
$$
Therefore $\linq\{x_{1},x_{2},\ldots\}$  is isomorphic to a Hilbert
space.
\eop
\bigskip

The equivalence of (i) and (ii) in the above corollary was first
proved by Arazy and Lindenstrauss \cite{AraLin}.
%\bigskip
%
%
% Section 7
%
%
\mysec{7}{Open questions}%
Here we gather some problems suggested by our work.
%
%
%	Problem 7.1
%
%
\begin{problem} \em
The method of proof employed in Proposition~\ref{B4}
suggests considering a version of property~$(M)$ with
respect to convergence in measure, meaning that $\limsup \|f+g_{n}\|$
depends only on $\|f\|$ whenever $g_{n}\to 0$ in measure. 
It might be interesting to pursue this idea further.

Let us
indicate why such a notion might be useful. Suppose $X$ is a function
space on $[0,1]$ that can be renormed to have this property. Then
every subspace generated by a sequence of pairwise disjoint functions
has the usual property~$(M)$. Consequently, by \cite[Th.~4.3]{Kal-M},
if an order continuous nonatomic Banach lattice $Y$
different from $L_{2}$ embeds isomorphically into $X$, then the Haar
system in $Y$ cannot embed on disjoint functions in $X$. 
This is an important
property in the investigation of rearrangement invariant spaces,
see \cite{JMST}, \cite{Kal-Mem} and \cite[Sect.~2.e]{LiTz2}. 

By way of illustration we observe that this line of reasoning
together with \cite[Th.~7.5]{Kal-Mem} implies the following result. 
\end{problem}
%
%
%  Prop 7.2
%
%
\begin{prop}\label{7.2}
If $X\neq L_{2}$ is a separable rearrangement invariant 
space of functions on $[0,1]$ that can be renormed to have the above
property~$(M)$ with respect to convergence in measure
and if $Y$ is a rearrangement invariant space on $[0,1]$ isomorphic to
$X$, then $X=Y$, and the norms of $X$ and $Y$ are equivalent (for
short, the space $X$ has unique rearrangement invariant
structure).
\end{prop}

One can check that the Lorentz spaces  $L_{p,q}[0,1]$ can be so
renormed and thus Proposition~\ref{7.2} can be applied to them.
In fact, the canonical norm works for $1\le q<p<\infty$, but for
$1< p<q<\infty$ the canonical expression is only a quasinorm,
and one should consider the norm
$$
\|f\|_{(p,q)} =
\left(\int_{0}^{1}\left(t^{1/p}f^{**}(t)\right)^{q}\frac{dt}{t}\right)
^{1/q}.
$$
Let us remark that $L_{p,q}[0,1]$ even has $(m_{q}(\tau))$  for
the topology $\tau$ of convergence in measure. For $q=\infty$ one
obtains this conclusion for the span of the bounded functions in
$L_{p,\infty}[0,1]$; note that this span is an $M$-ideal in its bidual
\cite{Dirk5}.

We conjecture that Orlicz spaces allow such a renorming, too. Anyway,
in the case of Orlicz spaces one can argue that a sequence of
pairwise disjoint functions generates a modular sequence space which can be
renormed to have property~$(M)$ by \cite[Prop.~4.1]{Kal-M}. 
This argument is a variant of that given in \cite[p.~168]{JMST} for the
reflexive case.
%
%
%	Problem 7.3
%
%
\begin{problem}	\em
Is $\K(\lp,c_p)$ an $M$-ideal in $\LL(\lp,c_p)$ for $1<p<2$ 
(cf.\ Proposition~\ref{B5}(a))? Note that $c_p$ is isomorphic to $\it
UT_p$ \cite{AraLin}.
\end{problem}
%
%
%	Problem 7.4
%
%
\begin{problem} \em
What about the commutative version of Proposition~\ref{B5}(a), i.e.,
is $\K(\lp,H_p)$ an $M$-ideal in $\LL(\lp,H_p)$ for $1<p<2$?
Recall that the Hardy space
$H_p$ must not be replaced by $\Lp$, by the results of
\cite{BangOd}.
\end{problem}
%
%
%	Problem 7.5
%
%
\begin{problem}	\em
By Corollary~\ref{G9} and well-known results on $M$-ideals (see
\cite{Kal-M} or \cite[Th.~V.5.4]{HWW})
there is a sequence $(K_n)$ of compact operators on the 
little Bloch space $\beta_0$ such that $K_n \to \Id$ strongly and
$\|\Id-2K_n\|\to 1$. Such operators can actually be obtained by
taking suitable convex combinations of Fej\'er operators. Is it
possible to find such $K_n$ explicitly among the classical kernel
operators? This would give a direct proof, which does not rely on the
embedding of $\beta_0$ into $c_0$, that $\K(\beta_0)$ is an
$M$-ideal in $\LL(\beta_0)$ and, moreover, that $\K(X,\beta_0)$ is an
$M$-ideal in $\LL(X,\beta_0)$ for all Banach spaces $X$, cf.\
\cite{AETD}. In fact, it suffices to find $K_n$ satisfying
$\|\Id-K_n\|\to 1$ \cite{Oja4}.   Also, it would be interesting to
check property~$(M^*)$ for $\beta_0$ directly.
\end{problem}
\bigskip

%
%
%  References
%
%
\typeout{References}
%\bibliography{basisref,dirkref,peterref,partbib}
%\bibliographystyle{standard}

%
%
% Address
%
%
\normalsize
\bigskip\bigskip\bigskip
\begin{tabbing}
\hspace*{.55 \textwidth} \= \kill
Department of Mathematics\> I.~Mathematisches Institut\\
University of Missouri  \> Freie Universit\"at Berlin \\
  \> Arnimallee 2--6 \\
Columbia, Mo 65\,211 \> 14\,195 Berlin \\[2pt]
U.S.A. \> Federal Republic of Germany			\\[4pt]
e-mail:\> e-mail: \\
mathnjk@mizzou1.missouri.edu  \>
werner@math.fu-berlin.de
\end{tabbing}

\end{document}